\documentclass[a4paper,11pt,twoside]{article}

\usepackage{xypic}
\usepackage{latexsym}
\usepackage{amsmath}
\usepackage{amstext}
\usepackage{amsfonts}
\usepackage{amssymb}
\usepackage{amsbsy}
\usepackage{eucal}
\usepackage{enumerate}
\usepackage{theorem}

\xyoption{all}
\CompileMatrices

\newtheorem{lemma}{Lemma}[section]
\newtheorem{proposition}[lemma]{Proposition}

\newtheorem{theorem}{Theorem}

\theorembodyfont{\rmfamily}
\newtheorem{remark}[lemma]{Remark}
\newtheorem{definition}[lemma]{Definition}
\newtheorem{example}[lemma]{Example}
\newtheorem{problem}[lemma]{Problem}

\newtheorem{assumption}[lemma]{Assumption}
\newtheorem{notation}[lemma]{Notation}

\newcommand{\id}{\mathop{\rm id}\nolimits}
\newcommand{\im}{\mathop{\rm Im}\nolimits}

\newcommand{\Spec}{\mathop{\rm Spec}\nolimits}

\newcommand{\Iso}{\mathop{\rm Iso}\nolimits}
\newcommand{\Gal}{\mathop{\rm Gal}\nolimits}

\newcommand{\Hom}{\mathop{\rm Hom}\nolimits}
\newcommand{\Pic}{\mathop{\rm Pic}\nolimits}

\newcommand{\Graph}{\mathop{\rm Graph}\nolimits}

\newcommand{\chr}{\mathop{\rm char}\nolimits}

\newcommand{\qed} {\hbox{} \nolinebreak \hfill $\;\Box$}

\begin{document}

\addtolength{\baselineskip}{0.5 pt}

\setlength{\parskip}{0.2 ex}

\title{Non-constant Curves of Genus 2 with Infinite Pro-Galois Covers}

\author{Claus Diem and Gerhard Frey}

\date{November 14, 2006}

\maketitle

\begin{abstract}
For every odd prime number $p$, we give examples of non-constant
smooth families of curves of genus $2$ over fields of
characteristic $p$  which have pro-Galois (pro-\'etale) covers of
infinite degree with geometrically connected fibers. The Jacobians
of the curves are isomorphic to products of elliptic curves.
\end{abstract}

\vspace{1 ex}
\begin{center}\begin{minipage}{110mm}\footnotesize{\bf Key words:} Covers of curves, fundamental groups, families of curves, curves of genus 2, curves with split Jacobian.
\end{minipage}
\end{center}

\begin{center}\begin{minipage}{110mm}\footnotesize{\bf MSC2000:} 14H30, 14H45, 14H10.

\end{minipage}
\end{center}
\vspace{1 ex}

\pagestyle{myheadings} \markboth{\sc Diem, Frey}
{\sc Curves of Genus 2 with Infinite Pro-Galois Covers}
\section{Introduction and Main Result}\label{intro}
In the following, unless stated otherwise, a \emph{curve} over a
field is assumed to be smooth, proper and geometrically connected.
This work is motivated by the following general problem.

Consider all curves of a
fixed genus $g$ over fields of a certain type (e.g.\ algebraically closed, local,
finitely generated). Among these curves, does there exist
a curve
$C$ over a field $K$ which allows an infinite tower of non-trivial
unramified covers $\cdots \longrightarrow C_i \longrightarrow \cdots \longrightarrow C_0 = C$ such that for all $i$, $C_i$ is
also a (\emph{geometrically connected}) curve over $K$ and $C_i
\longrightarrow C$ is Galois?

We concentrate on this problem for \emph{curves of genus $2$ over
finitely generated fields}. Some examples of curves of genus $2$
over \emph{finite fields} allowing such a tower of covers are
known, and we address the question whether there also exist
\emph{non-constant} curves of genus $2$ (over function fields)
with such a tower. We show with explicit examples that this is the
case.

\begin{example}\label{example1}
For every odd prime $p$, the (smooth, proper) curve of genus 2 given by the (affine) equation
\[ Y^2 = (t-1) \cdot (1-t^{p-1}) \cdot (X^2 -1) \cdot (X^2-t^{p-1}) \cdot (X^2 - 1 - t - \cdots - t^{p-1}) \]
over the global function field $K=\mathbb{F}_{p^2}(t)[\sqrt{t}, \sqrt{t-1}, \ldots, \sqrt{t-(p-1)}]$ allows
such a tower of covers.
\end{example}

\subsection{Further Background Information}
We will use the following terminology. A \emph{curve cover} over
$K$ is a surjective morphism of curves over $K$. A projective
limit of a projective system of curves (where the morphisms are curve covers) over $K$ is
called a \emph{pro-curve} over $K$. A \emph{pro-curve cover} of
a curve $C$ over $K$ is a surjective morphism $\pi : D
\longrightarrow C$ where $D$ is a pro-curve.
\begin{definition}
A curve cover is \emph{Galois} if it is unramified and the corresponding extension of function fields is Galois. A pro-curve cover is \emph{pro-Galois} if it is a projective limit of a projective system of Galois covers.
\end{definition}
For brevity, we speak of \emph{pro-Galois curve covers} instead of pro-Galois pro-curve covers.

With this terminology, the general problem posed above can be reformulated as follows:
\emph{Are there curves of a fixed genus $g$ over fields of certain types which allow pro-Galois curve covers of infinite degree?}

As an example consider the case that $K$ is an algebraically closed field. Then a curve $C$ over $K$
allows a pro-Galois curve cover of infinite degree and only if the genus of $C$ is \nolinebreak $\geq 1$.
 The situation is very different over finitely generated fields.

No curve of genus 1 over a finitely generated
field $K$ has a pro-Galois curve cover of infinite degree. There
even exists a universal bound $n=n(K)$ such that all Galois curve
covers of a genus 1 curve with a rational point over $K$ have degree $\leq n$.
Indeed, by the Theorem of Mazur-Kamienny-Merel (see \cite{Me}) and an
induction argument on the absolute transcendence degree of $K$, there exists a number
$n=n(K)$ such that all elliptic curves over $K$ have at most $n$ $K$-rational torsion
points -- from this the assertion follows; see \cite[Section 2]{FKV} for details.

If $K$ is finitely generated and $C$ is a curve over $K$ corresponding to the generic point of the
moduli scheme of curves of a certain genus, the authors expect that $C$ also does not have a
pro-curve cover of infinite degree. But for special
curves the situation changes. For every natural number $g \geq 3$ and every
field $K$ containing the 4-th roots of unity, explicit examples of curves over $K$ of
genus $g$ having a pro-Galois curve cover of infinite
degree can be given. If $K$ is not finite over its prime field the curves
can be chosen to be non-constant and even non-isotrivial
(i.e.\ for every extension of the ground field, they stay non-constant); see \cite{FKV, FK-proj}.

For finite fields, examples of curves of genus 2 with a pro-Galois
curve cover of infinite degree are known (see \cite[\S 3, Examples
2,3]{Ih}, \cite{FKV}), but to the knowledge of the authors, no
examples of non-constant curves of genus $2$ with a pro-Galois
curve cover have been known so far.

\subsection{The Main Result}
\label{subsection-main-result}

In this subsection we shall formulate the main result of this
paper. We thereby freely use some results of \cite{Di} which we recall at the beginning of Section \nolinebreak \ref{over Y_0(N)}. We begin with some notation.

Let $p$ be an odd prime, and let $N$ be an odd natural number.
Let $S$ be a smooth variety over a finite field of characteristic $p$. We assume that there is an $S$-isogeny $$\tau : \mathcal{E}
\longrightarrow \mathcal{E}'$$ between two elliptic curves over
$S$ having degree $N$. The kernel of the multiplication
by $2$ of $\mathcal{E}$ is denoted by $\mathcal{E}[2]$, the image of the zero-section by $[0_{\mathcal{E}}]$. We denote by $\mathcal{E}[2]^{\#}$ the
$S$-scheme $\mathcal{E}[2] - [0_{\mathcal{E}}]$. For the curve
$\mathcal{E}'$ we use similar definitions.

\begin{notation}\label{normalcurve}
Let  $\mathcal{P} := \mathcal{E}/\langle -1 \rangle,
\mathcal{P}' := \mathcal{E}'/\langle -1 \rangle$ with projections
$$\rho : \mathcal{E} \longrightarrow \mathcal{P} \,\,\mbox{and}\,\,\rho':
\mathcal{E}' \longrightarrow \mathcal{P}'.$$
\end{notation}
\begin{remark} The quotients $\mathcal{P}$ and $\mathcal{P}'$ are $\mathbb{P}^1$-bundles over $S$ which are isomorphic to $\mathbb{P}^1_S$ if $\mathcal{E}[2]$ is $S$-isomorphic
to the group scheme $(\mathbb{Z}/2\mathbb{Z})^2$. Without
weakening the results of the paper very much the reader may assume
this.
\end{remark}

There exists a unique $S$-isomorphism $$\gamma :
\mathcal{P} \tilde{\longrightarrow} \mathcal{P}'$$ such that
$$\rho' \circ \tau|_{\mathcal{E}[2]^{\#}} = \gamma \circ
\rho|_{\mathcal{E}[2]^{\#}}.$$

\begin{assumption}
Let the following two equivalent conditions be satisfied.
\label{assumption}
\begin{itemize}
\item
For no geometric point $s$ of $S$, there exists an
isomorphism $\alpha : \mathcal{E}_s \longrightarrow
\mathcal{E}_s'$ such that $\alpha|_{\mathcal{E}_s[2]} =
\tau_s|_{\mathcal{E}_s[2]}$.
\item
$\gamma(\rho([0_\mathcal{E}]))$ and $\rho'([0_{\mathcal{E}'}])$ are disjoint.
\end{itemize}
\end{assumption}

 \begin{notation}\label{curve}
Let $\mathcal{C}$ be the normalization of the integral scheme $\mathcal{E}
\times_{\mathcal{P}'} \mathcal{E}'$ (where the product is with
respect to $\gamma \circ \rho$ and $\rho'$).
\end{notation}
We thus have the commutative diagram
\[ \xymatrix{
 & \ar[dl] {\mathcal{C}} \ar[dr] & \\
\ar_{\rho}[d]  {\mathcal{E}} & & \ar^{\rho'}[d] {\mathcal{E}'} \\
\ar^{\gamma}[rr] {\mathcal{P}} & & {\mathcal{P}'\;\;.} }
\]

Now $\mathcal{C}$ is a curve of genus $2$ over $S$ and the morphisms $\mathcal{C} \longrightarrow \mathcal{E}$, $\mathcal{C} \longrightarrow \mathcal{E}$ are degree 2 covers. (For the terminology concerning relative curves, see Subsection \ref{terminology}.) It is the curve $\mathcal{C}$ which is the basic object of this article.

Our main result is the following theorem.

\begin{theorem}\label{main}
We use the notations from above, and we assume that Assumption \ref{assumption} is satisfied. So $\mathcal{C}$ is a curve of genus 2 over $S$ which is defined as in Notation \ref{curve}.

Then there exists
a connected Galois cover $T \longrightarrow S$ with Galois group a
(finite) elementary abelian $2$-group such that the curve
$\mathcal{C}_T$ over $T$ has a pro-Galois curve cover whose Galois
group $\mathcal{G}$ fits into an exact sequence
$$ 1 \longrightarrow \prod_{i=1}^\infty A_{n_i} \longrightarrow
\mathcal{G} \longrightarrow (\mathbb{Z}/2\mathbb{Z})^r
\longrightarrow 1 $$ with $r \leq 4$, $n_i$ are pairwise distinct
natural numbers and $A_{n_i}$ is the alternating group on $n_i$
elements.
\end{theorem}

\begin{remark}
It is easy to see that the Jacobian variety $J_{\mathcal{C}}$ of $\mathcal{C}$ is isogenous to $\mathcal{E} \times \mathcal{E}'$. We shall prove the stronger result that $J_{\mathcal{C}}$ is in fact \emph{isomorphic} to $\mathcal{E} \times \mathcal{E}'$; see Proposition \ref{isomorphic} as well as Subsection \ref{isogenous-isomorphic} for additional information.
\end{remark}

\subsection{Outline of the Proof}

The proof of Theorem \ref{main} is organized as follows:

In Section \ref{ramification},
 we show some basic results on ramification loci on
varying covers of a fixed elliptic curve by a fixed genus $2$
curve over an algebraically closed field. In Section \ref{branch},
we fix an isogeny $\tau : \mathcal{E}
\longrightarrow \mathcal{E}'$ as in Theorem \ref{main} (but over a more
general base scheme). Using results from \cite{Di}, we first show
that $\mathcal{C}$ is a curve of genus 2 whose Jacobian is
isomorphic to $\mathcal{E} \times \mathcal{E}'$. We then show with
the results of Section \ref{ramification} that there are
infinitely many covers from $\mathcal{C}$ to $\mathcal{E}$ with
the same ramification and branch loci and pairwise distinct
degrees. In Section \ref{Galois}, we show how our main result follows from the existence of infinitely many
minimal genus 2 covers of a given elliptic curve with the same
branch loci (under some conditions). 

Theorem \ref{main} follows by combining the final result of Section \ref{ramification}, Proposition \ref{sec3-result}, with the main result of Section \ref{Galois}, Proposition \ref{pro-Galois-general}.

\subsection{Derivation of Examples}
Theorem \ref{main} can be used to obtain Example \ref{example1}. Let $p$ be an odd prime, let $S :=
\mathbb{A}^1_{\mathbb{F}_p} - \{0 ,1, \ldots, p-1 \}$ and let us
denote the coordinate on $\mathbb{A}^1_{\mathbb{F}_p}$ by $t$. Let
$\mathcal{E}$ be the genus 1 curve over $S$ given by $$Y^2 Z -
X(X-Z)(X-tZ) = 0.$$ We want to fix a section of $\mathcal{E}$ over
$S$ in order to turn $\mathcal{E}$ into an elliptic curve. For
this, let $\mathcal{E}_a$ be the $S$-scheme given by $Y^2 -
X(X-1)(X-t) = 0$ and let $a : S \longrightarrow \mathcal{E}_a$ be
section given by $X \mapsto t, Y \mapsto 0$. We have a natural
inclusion $\iota: \mathcal{E}_a \hookrightarrow \mathcal{E}$ which
is compatible with the projection to $S$, and $\iota \circ a$ is a
section of $\mathcal{E}$ over $S$. We take this section as the
zero-section.

Similarly, let $\mathcal{E}'$ be the elliptic curve over $S$ given
by $$Y^2 Z - X(X-Z)(X-t^pZ) = 0$$ with zero-section determined by
$X \mapsto t^p, Y \mapsto 0$, and let
 $\tau : \mathcal{E} \longrightarrow \mathcal{E}'$ be the Frobenius endomorphism given by  $$X \mapsto X^p,\,\, Y \mapsto Y^p, Z \mapsto Z^p.$$

We note that the quotients $\mathcal{E} \longrightarrow
\mathcal{P}, \mathcal{E}' \longrightarrow \mathcal{P}'$ can be
identified with the usual projections to $\mathbb{P}^1_{S}$ (to
the ``$X$-coordinate''), and under these identifications,
$\gamma : \mathcal{P} \longrightarrow \mathcal{P}'$, becomes the identity. We also note that we have chosen
$S$ in such a way that the assumptions
of Theorem \ref{main} are satisfied. The curve $\mathcal{C}$ defined in Notation \ref{curve} is
the normalization of $\mathcal{E} \times_{\mathbb{P}^1_S}
\mathcal{E}'$. The field $K$ in the introduction is obtained as
the function field of the ``maximal'' connected Galois cover $T
\longrightarrow S$ which has an elementary abelian 2-group as
Galois group, that is, it is the maximal Galois extension of
$\mathbb{F}_p(t)$ which is unramified outside the places
corresponding to $0, 1, \dots, p-1$ and the place $\infty$ and
has an elementary abelian 2-group as Galois group. (The explicit description of the field can for example be obtained via Kummer theory.) The curve in
the introduction is the restriction of the curve $\mathcal{C}_T$
to the generic fiber of \nolinebreak $T$, that is, to $\Spec(K)$.
\medskip

Further examples can be obtained as follows:

Let $p$ be an odd prime and $N \geq 3$ and assume that $p$ and $N$ are coprime. It is well-known that there exists
a fine moduli
scheme for elliptic curves with cyclic subgroups of order $N$ and level-4-structure (of fixed determinant $\zeta_4$)
over schemes over
 $\mathbb{F}_{p}(\zeta_4)$. (We use that $p$ is odd.) Let us denote this scheme by
 $Y_0(N;4)_{\mathbb{F}_p(\zeta_4)}$ -- it is a smooth affine curve over $\mathbb{F}_p(\zeta_4)$.

Let $\tau : \mathcal{E} \longrightarrow \mathcal{E}'$ be the universal isogeny over
$Y_0(N;4)_{\mathbb{F}_p(\zeta_4)}$. There exists a uniquely determined open subscheme
$S$ of $Y_0(N;4)_{\mathbb{F}_p(\zeta_4)}$ such that the
assumptions of Theorem \ref{main} are satisfied. Obviously, $S$
is non-empty (as otherwise the universal elliptic curve over
$X_0(N;4)_{\mathbb{F}_p(\zeta_4)}$ would have complex
multiplication).

If we now apply Definition \ref{curve} to the restriction of the universal isogeny
$\tau : \mathcal{E} \longrightarrow \mathcal{E}'$ to $S$, we obtain a curve $\mathcal{C}$ over $S$. We can now apply Theorem \nolinebreak \ref{main} to this curve.

\subsection{Open Problems}

\begin{problem} It is an open problem whether there exists a curve
of genus 2 over a finitely generated field of characteristic 0
with a pro-Galois curve cover of infinite degree. The existence of
such a curve would of course be implied by the existence of a
genus 2 curve $\mathcal{C}$ with a pro-Galois curve cover of
infinite degree over an open part of $\text{Spec}(\mathcal{O}_K)$,
where $\mathcal{O}_K$ is the principal order in a number field. In
Subsection \ref{number-case} we discuss difficulties occurring
when one tries to adapt the proof of
Theorem \ref{main} to the
``mixed-characteristic case''.
\end{problem}

There is a sharpening of the condition on a curve over a field
having a pro-Galois curve cover of infinite degree.

Let $C$ be a curve over a field $K$ and let $P \in C(K)$. Following ideas of \cite{Ih},
the \emph{$K$-rational geometric fundamental group}
is introduced in \cite{FKV}. Let us denote this group by $\pi^{\text{\rm
geo}}_1(C,P)$.
The group $\pi^{\text{\rm geo}}_1(C,P)$ classifies pro-\'etale (curve) covers
 $c : D \longrightarrow C$ such that $P$ splits completely. If $\deg(c)$ is finite,
 this means that the topological point $P$ has exactly $\text{deg}(c)$ preimages under $c$.
  Thus the group $\pi^{\text{\rm geo}}_1(C,P)$ is infinite if and only if
  a pro-\'etale (or pro-Galois) (curve) cover of $D \longrightarrow C$ of infinite degree
  exists such that $P$ splits completely.

For every prime $p$, there exists a curve over a finite field $K$ of characteristic $p$ with
infinite $K$-rational geometric fundamental group. There also exists a non-isotrivial curve over
a function field in one variable over a finite field of given prime characteristic $p$ with infinite
$K$-rational fundamental group; this follows from \cite[Theorem 4.3]{FK-proj}.

Moreover, for every natural number $g \geq 3$ and every
field $K$ containing the 4-th roots of unity of characteristic congruent $3$ modulo $4$,
explicit examples of curves over $K$ of genus $g$ having an infinite $K$-rational fundamental
group can be given. If $K$ is not finite over its prime field the curves can be chosen to be non-isotrivial;
 see \cite[Theorem 1.1]{FK-proj}.

\begin{problem}
It is an open problem whether any curve over a finitely generated
field $K$ of characteristic 0 with infinite $K$-rational geometric
fundamental group exists,\footnote{The proof of \cite[Theorem
4.22]{FKV} and thus also the proof of the ``Result'' at the end of
the introduction of \cite{FKV} are incorrect: As stated in
\cite[Remark 5.4]{FK-proj}, the condition $\chr(K) \equiv 3 \text{
mod } 4$ has to be inserted.} and the case $g=2$ remains open in
any characteristic.
\end{problem}

\subsection{Terminology and Facts}
\label{terminology}
 The usage of ``Galois'' follows \cite[Expos\'e V]{SGA}. This means in
particular: If $Y$ is a locally noetherian, integral,
\emph{normal} scheme, $X$ is an integral scheme and $f: X
\longrightarrow Y$ is a finite surjective morphism, then $f$ is a
Galois cover if and only if it is unramified and the corresponding
extension of function fields is Galois; see \cite[Expos\'e I,
Corollaire 9.11]{SGA} and \cite[Expos\'e V, \S 8]{SGA}.

\smallskip

Let $S$ be a \emph{connected} locally noetherian scheme. A \emph{(relative) curve} of genus $g$ over $S$ is a smooth, proper $S$-scheme all of whose geometric fibers are connected non-singular curves. We denote curves over $S$ by $\mathcal{C}, \mathcal{D}, \mathcal{E}, \dots$. If $S$ is the spectrum of a field, we also use the notation $C, D, E, \dots$.

By the first corollary to the Theorem in \cite[Chapter II, \S 5]{Mu-AV}, the Euler-Poincar\'e characteristic of the geometric fibers of a relative curve is constant, thus the genus of the fibers is also constant; we define the \emph{genus} of a relative curve to be the genus of any of its fibers.

We will often use the fact that a curve over a regular scheme is regular; see
 \cite[Expos\'e II, Proposition 3.1]{SGA}.

A \emph{curve cover} of a curve $\mathcal{C}$ over $S$ is a finite and flat $S$-morphism
 $\pi: \mathcal{D} \longrightarrow \mathcal{C}$, where $\mathcal{D}$ is
 also a curve over $S$.

Note that an $S$-morphism $\pi : \mathcal{D} \longrightarrow \mathcal{C}$  between two curves over $S$ which induces curve covers on the fibers is automatically finite and flat, that is, it is a curve cover. (The morphism is fiberwise flat (\cite[Proposition 9.7]{Ha}) and thus flat by \cite[IV (11.3.11)]{EGA}, it is finite by \cite[Corollary 4.8]{Ha} and \cite[IV (8.11.1)]{EGA}.)

A \emph{pro-curve} over $S$ is a projective limit of a projective system of curves
with respect to morphisms which
are curve
covers over \nolinebreak $S$.

A \emph{pro-curve cover} of a curve $\mathcal{C}$ over $S$ is a projective limit of a projective system
of curve covers of $\mathcal{C}$. A pro-Galois pro-curve cover is abbreviated as
 \emph{pro-Galois curve cover}, and its automorphism group is called \emph{Galois group}.

\smallskip

We shall identify effective Cartier divisors on a locally noetherian integral scheme $X$ with the locally principal closed
subschemes of $X$; cf.\ \cite[Remark 6.17.1]{Ha}, \cite[21.2.12]{EGA}.
(This means in particular that Cartier divisors on an integral, regular schemes are identified with
closed subschemes of codimension 1.) If additionally
$X$ is an $S$-scheme, then under this identification, relative
effective Cartier divisors on $X$ over $S$ correspond to locally
principal closed subschemes of $X$ which are flat over $S$; see
\cite[Definition 3.4]{Mi-JV}.

\smallskip

Let $S$ be a scheme. If $S$ is integral, we denote its function field by $\kappa(S)$.
If $X$ is an $S$-scheme and $T \longrightarrow S$ is a morphism,
we denote the pull-back of $X$ to $T$ by $X_T$.
If $\alpha : X \longrightarrow Y$ is a morphism of $S$-schemes,
we denote its pull-back via $T$ by $\alpha_T : X_T \longrightarrow Y_T$.
Note that if $X$ is a curve over $S$ and $T$ is connected and locally noetherian, $X_T$ is a curve over $T$.

\smallskip

If $\mathcal{A}$ is an abelian scheme over $S$, the dual abelian scheme $\widehat{\mathcal{A}}$ exists by \cite[Theorem 1.9]{FC}. If $\varphi: \mathcal{A} \longrightarrow \mathcal{B}$ is a
homomorphism of abelian schemes over $S$, we denote the dual
homomorphism by $\widehat{\varphi} : \widehat{\mathcal{B}}
\longrightarrow \widehat{\mathcal{A}}$.

If $\mathcal{A}$ and $\mathcal{B}$
are two abelian schemes over $S$, we denote the
product over $S$ by
$\mathcal{A} \times \mathcal{B}$.
If $\mathcal{A}_1, \dots,
\mathcal{A}_n$ and $\mathcal{B}_1, \dots, \mathcal{B}_m$ are
abelian schemes over $S$, we have a canonical isomorphism between
the group of homomorphisms from
$\prod_j \mathcal{A}_j$
to
$\prod_i \mathcal{B}_i$ and the group of \emph{matrices}
$(a_{i,j})_{i=1 \dots, m, j=1, \dots, n}$, where $a_{i,j} \in
\Hom_S(\mathcal{A}_j,\mathcal{B}_i)$. The composition of
homomorphisms (if applicable) corresponds to the multiplication of
the matrices. If we use matrices, we identify $\mathbb{Z}$ with
its image in any endomorphism ring.

\smallskip

If $\mathcal{C}$ is a curve over $S$, we denote its Jacobian by
$J_{\mathcal{C}}$. The dual of the Jacobian is denoted by $J_{\mathcal{C}}$ and the canonical principal polarization on
$J_{\mathcal{C}}$ by $\lambda_{\mathcal{C}} : J_{\mathcal{C}}
\longrightarrow \widehat{J_{\mathcal{C}}}$; cf.\ \cite[Proposition 6.9]{Mu-GIT}. If $\pi : \mathcal{D} \longrightarrow \mathcal{C}$ is a curve cover over $S$, we have
canonical homomorphisms $\pi^* : J_{\mathcal{C}} \longrightarrow
J_{\mathcal{D}}$ and $\pi_* := \lambda_{\mathcal{C}}^{-1} \,
\widehat{\pi^*} \,
 \lambda_{\mathcal{D}} : J_{\mathcal{D}} \longrightarrow J_{\mathcal{C}}$.

\smallskip

Let $\mathcal{E}$ be an elliptic curve over $S$, $\mathcal{C}$ a
genus $2$ curve over $S$, $\pi : \mathcal{C} \longrightarrow
\mathcal{E}$ a curve cover. We call $\pi$ \emph{minimal} if it
does not factor through a non-trivial isogeny $\mathcal{E}_1
\longrightarrow \mathcal{E}$. Note that this is equivalent to
$\ker(\pi_* : J_{\mathcal{C}} \longrightarrow \mathcal{E})$ being
an elliptic curve, and it is equivalent to $\pi^* : \mathcal{E}
\longrightarrow J_{\mathcal{C}}$ being a closed immersion; see the
beginning of Section 2 of \cite{KaHur} and Point 7) in Section 7
of \cite{KaHur}.

\smallskip

Let $\pi : X \longrightarrow Y$ be a finite morphism of locally noetherian schemes.
Then the image $\pi(x) \in Y$ of a ramification point $x \in X$ is called \emph{branch point}.
The set of ramification points as well as the set of branch points are closed in $X$ and $Y$
respectively (the set of ramification points is closed because it the support of $\Omega_{X/Y}$, and the set of branch points is then the image of a closed set under a finite morphism, thus also closed). The corresponding schemes with the reduced induced scheme structures are called \emph{ramification locus} and \emph{branch locus} respectively.

\smallskip

We say that a field extension $K|k$ is \emph{regular}
if $K \otimes_k \overline{k}$ is a domain ($\overline{k}$ =
algebraic closure of $k$) (cf.\cite[VIII, \S 4]{La}). This notion
should not be confused with the notion of a regular scheme.

\subsection*{Acknowledgments}

The authors would like to thank E.\ Kani and E.\ Viehweg for
fruitful discussions on questions related to this work. They thank the anonymous referee for carefully reading the manuscript and for various helpful suggestions.

\section{The Ramification Loci of Covers of Genus $2$ of an Elliptic Curve}
\label{ramification}

Let $\overline{\kappa}$ be an algebraically closed field. When we
speak of the ``intersection'' of two curves on a smooth surface,
we mean the scheme-theoretic intersection. We will use some easy
results from intersection theory which we recall in an appendix to
this section.

Let $C$ be a curve of genus 2 over $\overline{\kappa}$.
\begin{proposition}
\label{unr-trans}
Let $P$ be a closed point of $C$, let $c : C \longrightarrow E$ be a cover which maps $P$ to $0_E$,
let $\iota : C \longrightarrow J_C$ be the canonical immersion which maps $P$ to $0_E$. Then the following assertions are equivalent.
\begin{itemize}
\item
$c$ is unramified in $P$.
\item
$\iota(C)$ and $\ker(c_*)$ intersect transversely in $\iota(P)$ (inside $J_C$).
\end{itemize}
\end{proposition}
\emph{Proof.} Let $c^{-1}(0_E) = C_{0_E} = C \times_E 0_E$ be the
fiber of $c$ at $0_E$. Note that $c = c_* \circ \iota$. Consider
the following Cartesian diagram.
\[\xymatrix{ c^{-1}(0_E) = C \times_E 0_E \ar@{^{(}->}[r] \ar@{->}[d] & C \ar@{^{(}->}^{\iota}[d] \\
\ker(c_*) \ar@{->}[d] \ar@{^{(}->}[r] & \ar@{->}[d]^{c_*}  J_C \\
0_E \ar@{^{(}->}[r] & E
}\]
We know that $c^{-1}(0_E)$ is a 0-dimensional closed subscheme on $C$, and its support contains $P$.
Both statements of the proposition are equivalent to the multiplicity of $c^{-1}(0_E)$ at $P$ being 1.
\qed

\begin{proposition}
\label{ramification-alternative}
Let $c_1 : C \longrightarrow E_1$, $c_2 : C \longrightarrow E_2$ be minimal covers such that there does
not
exist
an isomorphism of curves $\sigma : E_1
\tilde{\longrightarrow} E_2$ with $\sigma \circ c_{1} = c_{2}$.
Then $\ker(c_{1*})$ and $\ker(c_{2*})$ are abelian subvarieties of
$J_C$ which intersect in a finite group scheme. There is the
following alternative:
\begin{itemize}
\item
\emph{Either} the intersection of $\ker(c_{1*})$ and $\ker(c_{2*})$ is a reduced (\'{e}tale) finite group scheme
and the ramification loci of $c_1$ and $c_2$ are disjoint.
\item
\emph{or} the intersection of $\ker(c_{1*})$ and $\ker(c_{2*})$ is a non-reduced finite group scheme and
the ramification loci of $c_1$ and $c_2$ are equal.
\end{itemize}
Note that if $\chr(\overline{\kappa}) = 0$ the first alternative holds.
\end{proposition}
\emph{Proof.}
The fact that the abelian subvarieties $\ker(c_{1*})$ and $\ker(c_{2*})$ of $J_C$ intersect in a finite group scheme
follows because by assumption there does not exist a $\tilde{\sigma} \in \Iso(E_1,E_2)$ with
$\tilde{\sigma} \circ {c_1}_* = {c_2}_*$.

Note that $\ker(c_{1*})$ and $\ker(c_{2*})$ intersect transversely at 0 if and only if the intersection of the two
elliptic curves is a reduced (\'etale) finite group scheme; see also Lemma \ref{transversely-equiv}).

Let $P$ be a closed point of $C$. By composing the $c_i$ with suitable translations on $E_i$, we can assume that
$c_i(P) = 0_{E_i}$ $(i=1,2)$. We can now apply the previous proposition. By the fact that ``non-transversal intersection
at $P$'' is an equivalence relation (see Lemma \ref{non-transversality}), we conclude:

If $\ker(c_{1*}) \cap \ker(c_{2*})$ is a non-reduced finite group scheme, $P$ is a ramification point of $c_1$ if and
only if it is a ramification point of $c_2$.

If on the other hand $\ker(c_{1*}) \cap \ker(c_{2*})$ is a reduced finite group scheme, $P$ cannot be a ramification
point of both $c_1$ and $c_2$.
\qed

\begin{proposition}
\label{mod-p->ram-equality}
Assume that $\chr({\overline{\kappa}}) = p > 0$. Let $c_1 : C \longrightarrow E$ and $c_2 : C \longrightarrow E$ be
minimal covers. Let ${\overline{c}_1}_*, {\overline{c}_2}_* \in \Hom_{\overline{\kappa}}(J_C,E) \otimes_{\mathbb{Z}}
 \mathbb{F}_p$ be the induced elements, and assume that ${\overline{c}_1}_* = {\overline{c}_2}_*$. Then
 the ramification loci of $c_1$ and $c_2$ are equal.
\end{proposition}
\emph{Proof.} Under the assumption of the proposition, $\ker({c_1}_*)[p]$ is equal to $\ker({c_{2*}})[p]$
(as closed subschemes of $J_C$). In particular,  $\ker({{c_1}_*}) = \ker({{c_2}_*})$ or $\ker({{c_1}_*})
\cap \ker({{c_2}_*})$ is non-reduced. The result follows with the last proposition.
\qed

\smallskip

\subsection*{Appendix to Section \ref{ramification}: Facts about Intersection Theory}
Let $\overline{\kappa}$ be an algebraically closed field, let $X$ be a smooth surface over
${\overline{\kappa}}$, and let $P \in X$ be a closed point.
Let $Y$ be a curve  in $X$. As always we assume
that $Y$ is smooth.
For the following arguments it would suffice
that
$X$ and $Y$ are smooth at $P$.
Note that $X$ is locally factorial, so if $Y \hookrightarrow X$ is a 1-dimensional closed subscheme of $X$
containing $P$ with canonical closed immersion $\iota_Y : Y \hookrightarrow X$, the kernel of
$\iota_Y^{\#} : \mathcal{O}_{X,P} \longrightarrow \mathcal{O}_{Y,P}$ is generated by a single element
$f \in \mathcal{O}_{X,P}$ (unique up to a unit). We call such an $f$ a \emph{local equation} of $Y$ at $P$.
Note that $f \in \frak{m}_{X,P}$.
\begin{definition}
Let $Y_1$ and $Y_2$ be two curves on $X$ such that $P \in Y_1, P \in Y_2$. Then $Y_1$ and $Y_2$ \emph{intersect transversely at $P$} if the local
equations of $Y_1$ and $Y_2$ at $P$ generate $\frak{m}_{X,P}$.
\end{definition}
\begin{lemma}
\label{smooth-local-equation}
Let $Y$ be a curve in $X$, let $P$ be a closed point on $Y$. Then the local equation of $Y$ at $P$ does not lie in
$\frak{m}_{X,P}^2$.
\end{lemma}
\emph{Proof.}
This follows easily from the fact that $X$ is regular and 2-dimensional whereas $Y$ is regular and 1-dimensional. \qed

\smallskip

For the following lemma, note that the surjection $\frak{m}_{X,P} \longrightarrow \frak{m}_{Y_i,P}$ induces by dualization injections of (Zariski) tangent spaces $T_{Y_i,P} \hookrightarrow T_{X,P}$.

\begin{lemma}
\label{transversely-equiv}
Let $Y_i \subset X$ ($i=1,2$) be two curves in $X$ meeting in $P$. Let $Y_1 \cap Y_2 := Y_1 \times_X Y_2$ be the
scheme-theoretic intersection of $Y_1$ and $Y_2$. The following statements are equivalent:
\begin{enumerate}[a)]
\item
$Y_1$ and $Y_2$ intersect transversely at $P$.
\item
Viewed as elements in $\frak{m}_{X,P}/\frak{m}^2_{X,P}$, the local equations of $Y_1$ and $Y_2$ at $P$ are
linearly independent.
\item
The canonical homomorphism $\frak{m}_{X,P}/\frak{m}^2_{X,P} \longrightarrow \frak{m}_{Y_1,P}/\frak{m}^2_{Y_1,P} \times \frak{m}_{Y_2,P}/\frak{m}^2_{Y_2,P}$ is an isomorphism.
\item
The canonical homomorphism $T_{Y_1,P} \times T_{Y_2,P} \longrightarrow T_{X,P}$ is an isomorphism.
\item
There exists a neighborhood $U$ of $P$ such that $(Y_1 \cap Y_2)|_U$ is equal to the closed subscheme given by
the closed immersion $\Spec{({\overline{\kappa}})} \hookrightarrow U$ at $P$.
\end{enumerate}
\end{lemma}
\emph{Proof.}
a) and b) are equivalent by Nakayama's Lemma. The equivalence of b) and c) is easy, and d) is a ``dual formulation'' of c). The local ring of $Y_1 \cap Y_2$ at $P$ is isomorphic to $\mathcal{O}_{X,P}/(f_1, f_2)$. This implies that a) holds if and only if the local ring of $Y_1 \cap Y_2$ at $P$ is
isomorphic to $\Spec(\overline{\kappa})$. This in turn is equivalent to e).
\qed

\medskip

We define a relation on the set of all curves $Y$ lying in $X$ and meeting $P$ as follows: $Y_1$ is equivalent to
$Y_2$ if and only if $Y_1$ and $Y_2$ do not intersect transversely at $P$.

\begin{lemma}
\label{non-transversality}
The relation just defined is an equivalence relation.
\end{lemma}
\emph{Proof.}
We only have to prove the transitivity. This follows from point b) in the above equivalences with
Lemma \ref{smooth-local-equation}.
\qed

\section{Minimal Covers with Prescribed Branch Loci}\label{branch}
\label{over Y_0(N)}

In this section, we first show how results of \cite{Di} imply the statements which were implicitly used in Subsection \ref{subsection-main-result}. Then we use the results of the previous section to derive a preliminary result on our way to prove Theorem \nolinebreak \ref{main}.

Let us start off as in Subsection \ref{subsection-main-result}: Let $p$ be an odd prime number, and let $N$ be an odd natural
number. Let $S$ be a locally noetherian, integral, regular scheme of characteristic
$p$. Let $\tau : \mathcal{E} \longrightarrow \mathcal{E}'$ be an
isogeny of elliptic curves over $S$ of degree $N$. 

Also as in Subsection \ref{subsection-main-result}, let $\mathcal{E}[2]^{\#} := \mathcal{E}[2] -
[0_{\mathcal{E}}]$, where $[0_{\mathcal{E}}]$ is the Cartier
divisor associated to the zero-section $0_\mathcal{E}$ of
$\mathcal{E}$ (similar definitions for $\mathcal{E}'$), let
$\mathcal{P} := E/\langle -1 \rangle, \mathcal{P}' := E'/\langle
-1 \rangle$ with projections $\rho : \mathcal{E} \longrightarrow
\mathcal{P}$ and $\rho': \mathcal{E}' \longrightarrow
\mathcal{P}'$.

Let $\psi := \tau|_{\mathcal{E}[2]} : \mathcal{E}[2] \longrightarrow \mathcal{E}'[2]$.

It is
well-known
that $\mathcal{P}$ is isomorphic to the
$\mathbb{P}^1$-bundle
$\mathbb{P}(q_*(\mathcal{L}(0_{\mathcal{E}}))$, where $q :
\mathcal{E} \longrightarrow S$ is the structure morphism (see for example the discussion above Theorem 2 in \cite{Di}), and a similar statement holds for
$\mathcal{P}'$. Note that $\mathcal{E}[2]^{\#} \longrightarrow S$
and $\mathcal{E}'[2]^{\#} \longrightarrow S$ are \'etale covers of
degree 3 and $\mathcal{E}[2]^{\#} \longrightarrow \mathcal{P}$ as
well as $\mathcal{E}'[2]^{\#} \longrightarrow \mathcal{P}'$ are
closed immersions. By \cite[Propsition B.4]{Di}, there exists a
unique $S$-isomorphism $\gamma : \mathcal{P} \longrightarrow
\mathcal{P}'$ with $\gamma \circ \rho|_{\mathcal{E}[2]^{\#}} =
\rho' \circ \psi$. (The fact that $\mathcal{P}$ is isomorphic to
$\mathbb{P}^1_S$ if $\mathcal{E}[2] \simeq
(\mathbb{Z}/2\mathbb{Z})^2$ also follows from \cite[Propsition
B.4]{Di}.)

\begin{lemma}\label{empty} The conditions 
\begin{itemize}
\item
for no geometric point $s$ of $S$, there exists an isomorphism $\alpha : \mathcal{E}_s \longrightarrow \mathcal{E}_s'$
such that $\alpha|_{\mathcal{E}_s[2]^{\#}} = \psi$
\item
$\gamma \circ \rho([0_\mathcal{E}])$ and $\rho'([0_{\mathcal{E}'}])$ are disjoint
\end{itemize}
are equivalent.
\end{lemma}
\emph{Proof.} It is enough to prove the lemma in the case that $S$ is the spectrum of an algebraically closed field field. In this
situation one can use \cite[Lemma A.1]{Di}. \qed

\medskip
\emph{For the rest of this section, we assume that the conditions of Lemma \ref{empty} are satisfied.}
\begin{lemma}
\label{integral-lemma}
$\mathcal{E} \times_{\mathcal{P}'} \mathcal{E}'$ is integral.
\end{lemma}
As in Notation \ref{normalcurve},
 let $\mathcal{C}$ be the
normalization of $\mathcal{E} \times_{\mathcal{P}'} \mathcal{E}'$, and let $$\pi : \mathcal{C} \longrightarrow \mathcal{E}, \pi' : \mathcal{C} \longrightarrow \mathcal{E}'$$ be the canonical projections.

\begin{proposition}
\label{C-is-genus2}
$\mathcal{C}$ is a curve of genus 2 over $S$ and the maps $\pi, \pi'$ are degree 2 covers.
\end{proposition}
For Lemma \ref{integral-lemma} and Proposition \ref{C-is-genus2}, see \cite[Proposition 4.4]{Di}.

\medskip

In the notation of \cite{Di} (\cite[Definition 2.7, Proposition 2.15]{Di}), $(\mathcal{C},\pi,\pi')$
is a ``normalized symmetric pair'' corresponding to $(\mathcal{E}, \mathcal{E}', \psi)$. This means that
\begin{itemize}
\item
$\ker(\pi_*) = \im((\pi')^*)$ and $\ker(\pi^*) = \im(\pi'_*)$,
\item
$\pi_*W_{\mathcal{C}} = \mathcal{E}[2]^{\#}$ and $\pi'_*W_{\mathcal{C}} = \mathcal{E}'[2]^{\#}$,
\item
$\pi^*|_{\mathcal{E}[2]} = (\pi')^* \circ \psi$.
\end{itemize}
Here, $W_{\mathcal{C}}$ is the Weierstra{\ss} divisor of $\mathcal{C}$
and $\mathcal{E}[2]^{\#} := \mathcal{E}[2] - [0_{\mathcal{E}}]$, where $[0_{\mathcal{E}}]$ is the
Cartier divisor associated to the zero-section $0_{\mathcal{E}}$ of $\mathcal{E}$
(similar definition for $\mathcal{E}'[2]^{\#}$).

\begin{remark}
Even without the assumption that $S$ is regular, one can prove
that a genus 2 curve $\mathcal{C}$ over $S$ and two degree $2$
covers $\pi, \pi'$ with the above conditions exist. Moreover, the
triple $(\mathcal{C},\pi,\pi')$ is up to unique isomorphism (defined in an obvious way) uniquely determined by the three conditions. On the other hand, if the two equivalent
conditions in Lemma \ref{empty} are not satisfied, no
normalized symmetric pair corresponding to $(\mathcal{E},
\mathcal{E}, \psi)$ exists. For more information on these issues,
we refer the reader to \cite{Di}.
\end{remark}

Let $\delta :  \mathcal{E} \times \mathcal{E}' \longrightarrow
J_{\mathcal{C}}$ be the homomorphism given by the matrix $\left(
\begin{array}{cc} \pi^* & (\pi')^* \end{array} \right)$. (We
identify $\mathcal{E}$ and $\mathcal{E}'$ with their  dual abelian varieties.) 

By \cite[Proposition 2.14]{Di}, the kernel of $\delta$ is $\Graph(-\psi) = \Graph(\psi)$,
and the pull-back of the canonical principal polarization of $J_{\mathcal{C}}$ via $\pi$ is
twice the canonical product polarization of $\mathcal{E} \times \mathcal{E}'$.

\begin{proposition}
\label{isomorphic}
The Jacobian $J_{\mathcal{C}}$ is isomorphic to $\mathcal{E} \times \mathcal{E}'$.
\end{proposition}
\emph{Proof.} As said above, $J_{\mathcal{C}}$ is isomorphic to
$(\mathcal{E} \times \mathcal{E}')/\Graph(-\psi)$. The latter is
in turn isomorphic to $\mathcal{E} \times \mathcal{E}'$. In fact,
the isogeny $$\Phi : \mathcal{E} \times \mathcal{E}'
\longrightarrow \mathcal{E} \times \mathcal{E}'$$ given by the
matrix $\left( \begin{array}{cc} 2 & 0 \\ \tau & 1 \end{array}
\right)$ has kernel $\Graph(-\psi)$, thus it induces an
isomorphism $(\mathcal{E} \times \mathcal{E}')/\Graph(-\psi)
\tilde{\longrightarrow} \mathcal{E} \times \mathcal{E}'$.
\qed

\medskip

We use the group structure on $\mathcal{E}$
and obtain for all $a,b \in \mathbb{Z}$
  a morphism $a \pi + b
\widehat{\tau} \pi' : \mathcal{C} \longrightarrow \mathcal{E}$. We
have a bilinear form
\[ \beta: \; \mathbb{Z}^2 \times \mathbb{Z}^2 \longrightarrow \text{End}(\mathcal{E}), \;
((a,b),(c,d)) \mapsto (a \pi + b \widehat{\tau} \pi')_* \, (c {\pi} + d \widehat{\tau} {\pi'})^*\]
with $\beta((a,b),(a,b)) = \deg(a \pi + b \widehat{\tau} \pi')$.
As by definition $\pi_* \circ (\pi')^*  = 0$ and $\pi'_* \circ  \pi^* = 0$,
we have $\beta((a,b),(c,d)) = 2ac + 2Nbd \in \mathbb{Z}$. This implies that
\begin{equation}
\label{deg}
\deg(a\pi + b\widehat{\tau} \pi') = 2a^2 + 2Nb^2.
\end{equation}
Let us fix the following notation: For $(a,b) \in \mathbb{Z}^2$,
we denote the ramification locus of $a\pi + b \widehat{\tau} \pi'$
by $V_{(a,b)}$ and the branch locus by $\Delta_{(a,b)}$. Further,
we set $V := V_{(1,0)}$ and $\Delta := \Delta_{(1,0)}$. The
following lemma is \cite[Proposition 3.11]{Di}.

\begin{lemma}
\label{zero-element}
$\pi'|_V$ is the zero-element in the abelian group $\mathcal{E}'(V)$.
\end{lemma}

\begin{proposition}
\label{good-covers}
If $p | b$, then $V_{(1,b)} = V$ and $\Delta_{(1,b)} = \Delta$.
\end{proposition}
\emph{Proof.} Let $\overline{\kappa}$ be any algebraically closed field, and let $s$ be any point of $S(\overline{k})$. By Proposition \ref{mod-p->ram-equality},  $\pi_s$ and $(\pi + b \widehat{\tau} \pi')_s : \mathcal{C}_s \longrightarrow \mathcal{E}_s$ have the same ramification locus.

Now let $x$ be a $\overline{\kappa}$-valued point of $\mathcal{C}$ lying over $s$. Then the support of $x$ is a ramification point of $\pi$ if and only if it is a ramification point of $\pi_s$: both conditions are equivalent to the inclusion of $\overline{\kappa}$ into the local ring of $x$ in the geometric fiber $\mathcal{C}_{\pi(x)}$ being not surjective. A corresponding statement holds for $\pi + b \widehat{\tau} \pi'$.

This implies that $\pi$ and $\pi + b \widehat{\tau} \pi' : \mathcal{C} \longrightarrow \mathcal{E}$ have the same ramification locus, that is, the underlying sets of $V_{(1,b)}$ and $V$ are equal. As both schemes are endowed with the reduced induced scheme structure, we have that $V_{(1,b)} = V$.

By Lemma \ref{zero-element}, it follows that $$(\pi + b \widehat{\tau}
\pi')|_V = \pi|_{V} + b \widehat{\tau} \pi'|_V = \pi|_V,$$
 In
particular, $V$ is mapped under $\pi + b \widehat{\tau} \pi'$ to
$\Delta$. This implies that $\Delta_{(1,b)} = \Delta$. \qed

\begin{proposition}
If $2 | b$, then $\pi + b \widehat{\tau}  \pi'$ is minimal of degree $2 + 2Nb^2$.
\end{proposition}
\emph{Proof.}
The statement on the degree follows from (\ref{deg}).

To show that $\pi + b \widehat{\tau}  \pi'$ is minimal, that is, that
it does not factor through a non-trivial isogeny $\widetilde{\mathcal{E}}
\longrightarrow \mathcal{E}$, it suffices to show that the homomorphism
$$\pi_* + b \widehat{\tau}  \pi'_* : J_{\mathcal{C}} \longrightarrow \mathcal{E}$$
 does not factor through a non-trivial isogeny $\widetilde{\mathcal{E}} \longrightarrow \mathcal{E}$.

So assume that $\pi_* + b \widehat{\tau}$ factors through the isogeny $\widetilde{\mathcal{E}} \longrightarrow \mathcal{E}$. Then 
$$(\pi_* + b \widehat{\tau} \pi_*') \circ \pi^* : \mathcal{E} \longrightarrow \mathcal{E}$$ 
also
factors through $\widetilde{\mathcal{E}} \longrightarrow \mathcal{E}$. Now $(\pi_* + b \widehat{\tau} \pi_*') \circ \pi^* = 2 \cdot \id_{\mathcal{E}}$. This implies that the degree of the isogeny $\tilde{\mathcal{E}} \longrightarrow \mathcal{E}$ divides 4.

To rule out that the degree is 2 or 4, we consider the commutative diagram
\[ \xymatrix{
{\mathcal{E} \times \mathcal{E}'} \ar_{\delta}[d]  \ar^{\Phi}[dr] \ar@/_15ex/[dd]_{[2]}\\
{J_{\mathcal{C}} \simeq (\mathcal{E} \times \mathcal{E}')/\Graph(-\psi)} \ar[r]^>>>>>{\sim} \ar_{\widehat{\delta} \,
 \lambda_{\mathcal{C}}} [d] \ar@/_10ex/[dd]_{\pi_* +  b \widehat{\tau} \pi'_*} & {\mathcal{E} \times \mathcal{E}'\, ,}
 \ar[dl]^{\Psi}  \\
{\mathcal{E} \times \mathcal{E}'} \ar[d] \\
{\mathcal{E}}
}\]
where $\Phi$ is given as in Proposition \ref{isomorphic}, $\Psi : \mathcal{E} \times \mathcal{E}' \longrightarrow \mathcal{E} \times \mathcal{E}'$
is given by the matrix $\left( \begin{array}{cc} 1 & 0 \\ - \tau & 2 \end{array} \right)$, and the last
vertical arrow is given by the matrix $(\begin{array}{cc} 1 & b \widehat{\tau} \end{array})$. The homomorphism
$$\widehat{\delta} \, \lambda_{\mathcal{C}} : J_{\mathcal{C}} \longrightarrow \mathcal{E} \times \mathcal{E}'$$
 is given by the matrix $\left( \begin{array}{c} \pi_* \\ \pi'_* \end{array} \right)$.

Under the horizontal isomorphism in the diagram, $\pi_* + b \widehat{\tau} \pi'_*$ is given by the matrix $(\begin{array}{cc} 1 - bN & 2b \widehat{\tau} \end{array})$. Let $\iota : \mathcal{E} \longrightarrow \mathcal{E} \times \mathcal{E}' \simeq J_{\mathcal{C}}$ be the inclusion of the first summand. Then $( \pi_* + b \widehat{\tau} \pi'_*) \circ \iota = (1 - bN) \cdot \id_{\mathcal{E}}$. As by assumption $1 - bN$ is odd, the degree of
$\widetilde{\mathcal{E}} \longrightarrow \mathcal{E}$ cannot be
divisible by 2. It is thus an isomorphism.

We conclude that $\pi + b \widehat{\tau} \pi' : \mathcal{C} \longrightarrow \mathcal{E}$ is minimal. \qed

\medskip

The last two propositions imply:

\begin{proposition}
\label{result-for-section-1}
The covers $\pi_i := \pi + 2ip \widehat{\tau}  \pi' : \mathcal{C} \longrightarrow \mathcal{E}$ ($i \in \mathbb{N}$)
are minimal of degree $2 + 8N(ip)^2$ and have the same ramification loci and the same branch loci.
\end{proposition}

We summarize the results stated in Lemma \ref{empty},
  Lemma \ref{integral-lemma} and
in Propositions \ref{C-is-genus2}, \ref{isomorphic} and
\ref{result-for-section-1}
in the following proposition.
\begin{proposition}
\label{sec3-result}
Let $p$ be an odd prime, let $N$ be odd. Let $S$
be a locally noetherian, integral, regular scheme of
characteristic $p$, let $\tau : \mathcal{E} \longrightarrow
\mathcal{E}'$ be a cyclic isogeny of degree $N$ over $S$.

Then there exists a unique $S$-isomorphism $\gamma : \mathcal{P}
\tilde{\longrightarrow} \mathcal{P}'$ such that $\rho' \circ
\tau|_{\mathcal{E}[2]^{\#}} = \gamma \circ
\rho|_{\mathcal{E}[2]^{\#}}$. Assume that the following two
equivalent conditions are satisfied.
\begin{itemize}
\item
For no geometric point $s$ of $S$, there exists an isomorphism
$\alpha : \mathcal{E}_s \longrightarrow \mathcal{E}_s'$
 such that $\alpha|_{\mathcal{E}_s[2]} = \tau_s|_{\mathcal{E}_s[2]}$.
\item
$\gamma(\rho([0_\mathcal{E}]))$ and $\rho'([0_{\mathcal{E}'}])$
are disjoint.
\end{itemize}
 Let $\mathcal{C}$ be the normalization of the integral scheme
 $\mathcal{E} \times_{\mathcal{P}'} \mathcal{E}'$ (where the product is with respect to $\gamma \circ \rho$ and $\rho'$).
  Then $\mathcal{C}$ is a curve of genus 2 with $J_{\mathcal{C}} \simeq \mathcal{E} \times \mathcal{E}'$,
  the canonical morphisms $\pi : \mathcal{C} \longrightarrow \mathcal{E}$ and
  $\pi' : \mathcal{C} \longrightarrow \mathcal{E}'$ are degree 2 covers and there exists a
  sequence $(\pi_i)_{i \in \mathbb{N}_0}$ of minimal covers $\pi_i : \mathcal{C} \longrightarrow \mathcal{E}$
  with pairwise distinct degrees and $\pi_0 = \pi$ such that the ramification loci as well as the branch loci of the $\pi_i$ are all equal.
\end{proposition}

\section{Pro-Galois Curve Covers of Infinite Degree}\label{Galois}
\label{infinite-degree}
The goal of this section is to prove the following proposition
which together with Proposition \ref{sec3-result} implies  Theorem \ref{main}.

\begin{proposition}
\label{pro-Galois-general} 
Let $S$ be an integral, regular scheme of finite type over $\mathbb{Z}[1/2]$ or more generally a locally noetherian, integral, regular scheme over $\mathbb{Z}[1/2]$ such that $\pi_1(S)/\pi_1(S)^2$ is finite. (Here $\pi_1(S)$ denotes the fundamental group of $S$ with respect to some base point.)

Let $\mathcal{E}$ be an elliptic curve over $S$. Assume that there
exists a sequence $(\pi_i)_{i \in \mathbb{N}_0}$ of minimal covers
$\pi_i : \mathcal{C}_i \longrightarrow \mathcal{E}$ (where the
$\mathcal{C}_i$ are curves of genus $2$ over $S$) with pairwise
distinct degrees and $\deg(\pi_0) = 2$ as well as $\deg(\pi_i) \geq 5$ for $i \geq 1$ such that the branch loci of the $\pi_i$ are all equal.

Let $\mathcal{C} := \mathcal{C}_0$. Then there exists a connected Galois cover
$T \longrightarrow S$ with Galois group a (finite) elementary abelian $2$-group such
that the curve $\mathcal{C}_T$ over $T$ has a pro-Galois curve cover whose
Galois group $\mathcal{G}$ fits into an exact sequence
\[ 1 \longrightarrow \prod_{i=1}^\infty A_{n_i} \longrightarrow \mathcal{G}
\longrightarrow (\mathbb{Z}/2\mathbb{Z})^r \longrightarrow 1 \]
for some $r \leq 4$, where $n_i := \deg(\pi_i)$ and
$A_{n_i}$ is the alternating group on $n_i$ elements.
\end{proposition}

The rest of this section is devoted to a proof of this proposition.

\medskip

We first show that any integral, regular scheme $S$ of finite type over $\mathbb{Z}[1/2]$ the assumption that $\pi(S)/\pi(S)^2$ is finite is satisfied. We thereby give some background information on this condition as well.

Let $S$ be any scheme over $\mathbb{Z}[1/2]$. Then we have an exact sequence
\[ 0 \longrightarrow  \Gamma(S,\mathcal{O}_S)^*/\Gamma(S,\mathcal{O}_S)^{*2} \longrightarrow \Hom(\pi_1(S),\mathbb{Z}/2\mathbb{Z}) \longrightarrow \Pic(S)[2] \longrightarrow 0 \; ;\]
see \cite[Proposition 4.11]{Mi-EC} together with \cite[Corollary 4.7]{Mi-EC}. We thus see that $\pi_1(S)/\pi_1(S)^2$ is finite if and only if both $\Pic(S)[2]$ and $\Gamma(S,\mathcal{O}_S)^*/\Gamma(S,\mathcal{O}_S)^{*2}$ are finite. 

If now $S$ is an integral, normal scheme of finite type over $\mathbb{Z}$, both groups $\Pic(S)[2]$ and $\Gamma(S,\mathcal{O}_S)^*/\Gamma(S,\mathcal{O}_S)^{*2}$ are in fact finite, and hence so is $\pi_1(S)/\pi_1(S)^2$; see \cite[Theorems 7.4 and 7.5]{Lang-fundamentals}.

\medskip

Now let the assumptions of the proposition be satisfied. Let $\Delta$ be the branch locus of each of the $\pi_i$. As all
residue characteristics are distinct from 2 and $\deg(\pi_0) =2$,
for any $s \in S$ $(\pi_0)_s$ is generically \'etale, in
particular, $(\pi_0)_s$ is finite, and all branch points in
codimension 1 on $\mathcal{E}$ lie in the generic fiber over $S$.

More precisely, we have by \cite[Lemma 3.13]{Di}:

\begin{lemma}
The canonical morphism $\Delta \longrightarrow S$ is an \'etale cover of degree \nolinebreak 2.
\end{lemma}

Note that the curve $\mathcal{E}$ over $S$ gives rise to an
elliptic curve over $\kappa(S)$. In particular, the field extension $\kappa(\mathcal{E})|\kappa(S)$ is regular.
Let $L_i|\kappa(\mathcal{E})$ be the Galois
closure of the extension of function fields $\pi_i^{\#} : \kappa(\mathcal{E}) \hookrightarrow \kappa(\mathcal{C}_i)$.

Let us fix some compositum $L_i \overline{\kappa(S)}$ of $L_i$ and $\overline{\kappa(S)}$ over $\kappa(S)$.
\begin{lemma}
The Galois group of $L_i \overline{\kappa(S)}|\kappa(\mathcal{E}) \overline{\kappa(S)}$ is isomorphic
to the symmetric group $S_{n_i}$. In particular, the Galois group of $L_i|\kappa(\mathcal{E})$ is
isomorphic to $S_{n_i}$, and $L_i|\kappa(S)$ is regular.
\end{lemma}
\emph{Proof.}
As the extension of function fields $\kappa(\mathcal{C}_i) \overline{\kappa(S)}|\kappa(\mathcal{E}) \overline{\kappa(S)}$
 over $\overline{\kappa(S)}$ has 2 branched places, the conorm of each of these places has the form
 $\frak{P}_1^2 \frak{P}_2 \cdots \frak{P}_{n_i-1}$, where the $\frak{P}_j$ are pairwise distinct places. It follows that the Galois group of $L_i \overline{\kappa(S)}|\kappa(\mathcal{E})
\overline{\kappa(S)}$ (seen as permutation group on $n_i$ elements) contains a transposition.

Moreover, as $\pi_i : \mathcal{C}_i \longrightarrow \mathcal{E}$
is minimal, $\kappa(\mathcal{C}_i)|\kappa(\mathcal{E})$ has no
proper intermediate fields, and the same is true for
$\kappa(\mathcal{C}_i) \overline{\kappa(S)}|\kappa(\mathcal{E})
\overline{\kappa(S)}$.

Hence the Galois group of $L_i \overline{\kappa(S)}|\kappa(\mathcal{E}) \overline{\kappa(S)}$ is a primitive transitive subgroup of $S_{n_i}$ with a transposition, and hence it is equal to $S_{n_i}$; see \cite[Theorem 13.3]{Wi}. \qed

\medskip

Let $L_i^0 := L_i^{A_{n_i}}$. Let $\widetilde{\mathcal{C}}_i$ be
the normalization of $\mathcal{E}$ in $L_i$, and let
$\widetilde{\mathcal{C}}_i^0$ be the normalization of
$\mathcal{E}$ in $L_i^0$.

Let us fix inclusions of $\kappa(\widetilde{\mathcal{C}}_i)|\kappa(\mathcal{E})$ into some algebraic closure of $\kappa(\mathcal{E})$. For $t \in \mathbb{N}$, let $\mathcal{D}_t$ be the normalization
of $\mathcal{E}$ in the compositum of
$\kappa(\widetilde{\mathcal{C}}_0), \dots,
\kappa(\widetilde{\mathcal{C}}_t)$ over $\kappa(\mathcal{E})$, and
let $\mathcal{D}_t^0$ be the normalization of $\mathcal{E}$ in
the compositum of $\kappa(\widetilde{\mathcal{C}}_0^0), \dots,
\kappa(\widetilde{\mathcal{C}}_t^0)$ over $\kappa(\mathcal{E})$ (both composita with respect to the inclusions into the fixed algebraic closure of $\kappa(\mathcal{E})$). The extensions $\kappa(\mathcal{D}_t)|\kappa(\mathcal{E})$ and
$\kappa(\mathcal{D}_t^0)|\kappa(\mathcal{E})$ are Galois, and
$\kappa(\mathcal{D}_t^0)|\kappa(\mathcal{E})$ has as Galois group
an elementary abelian 2-group.
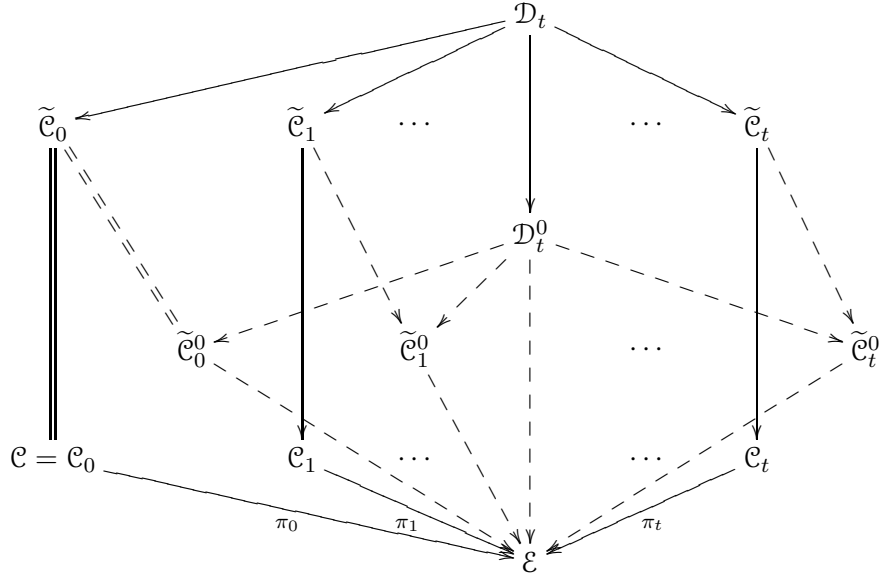
\[ \xymatrix{
& & & & {\mathcal{D}_t} \ar[lllld] \ar[lld] \ar[rrd] \ar[dd]\\
{\widetilde{\mathcal{C}}_0} \ar@{=}[ddd] \ar@{==}[ddr]& & {\widetilde{\mathcal{C}}_1} \ar[ddd] \ar@{-->}[ddr] & \cdots & & \cdots & {\widetilde{\mathcal{C}}_t} \ar[ddd] \ar@{-->}[ddr] \\
& & & & {\mathcal{D}^0_t} \ar@{-->}[ddd] \ar@{-->}[dlll] \ar@{-->}[dl] \ar@{-->}[drrr] \\
& {\widetilde{\mathcal{C}}_0^0} \ar@{-->}[rrrdd] & & {\widetilde{\mathcal{C}}_1^0} \ar@{-->}[rdd] & &  \cdots & &  {\widetilde{\mathcal{C}}_t^0} \ar@{-->}[ddlll] \\
{\mathcal{C} = \mathcal{C}_0} \ar[rrrrd]_{\pi_0} & & {\mathcal{C}_1} \ar[rrd]_{\pi_1} & \cdots & & \cdots & {\mathcal{C}_t} \ar[lld]^{\pi_t} \\
& & & & {\mathcal{E}}
}\]

Our goal is now to prove the following proposition.
\begin{proposition}\label{final}
\label{we-have-to-show}
For all $t \in \mathbb{N}$, we have:
\begin{enumerate}[a)]
\item
The morphisms $\mathcal{D}_t \longrightarrow \mathcal{C}$ and $\mathcal{D}_t^0 \longrightarrow \mathcal{C}$
are Galois covers.
\item
The Galois group of $\mathcal{D}_t^0 \longrightarrow \mathcal{C}$ is an elementary abelian 2-group.
\item
The restrictions $\Gal(\kappa(\mathcal{D}_t)|\kappa(\mathcal{D}^0_t)))
 \longrightarrow \Gal(\kappa(\widetilde{\mathcal{C}}_i)|\kappa(\widetilde{\mathcal{C}}_i^0))$
induce an isomorphism
\[\Gal(\mathcal{D}_t \rightarrow \mathcal{D}^0_t) \simeq \Gal(\kappa(\mathcal{D}_t)|\kappa(\mathcal{D}^0_t)) \tilde{\longrightarrow} \prod_{i=1}^{t} \Gal(\kappa(\widetilde{\mathcal{C}}_i)|\kappa(\widetilde{\mathcal{C}}_i^0)) \approx \prod_{i=1}^t A_{n_i}.\]
\item
Let $F_t$ be the algebraic closure of $\kappa(S)$ in $\kappa(\mathcal{D}^0_t)$, let $S_t^0$ be
the normalization of $S$ in $F_t$. Then $S_t^0 \longrightarrow S$ is a Galois cover, and $\mathcal{D}_t^0$ as well as $\mathcal{D}_t$ are in a canonical way curves over $S_t^0$. The induced morphism $\mathcal{D}_t^0 \longrightarrow \mathcal{C}_{S_t^0}$ is a Galois curve cover over $S_t^0$ with Galois group isomorphic to $(\mathbb{Z}/2\mathbb{Z})^r$ for some $r \leq 4$. We have an exact sequence of Galois groups
\[ 1 \longrightarrow \Gal(\mathcal{D}_t^0 \rightarrow \mathcal{C}_{S^0_t}) \longrightarrow \Gal(\mathcal{D}_t^0 \rightarrow \mathcal{C}) \longrightarrow \Gal(S_t^0 \rightarrow S) \longrightarrow 1. \]
In particular, $\Gal(S_t^0 \rightarrow S)$ is an elementary abelian 2-group.
\end{enumerate}
\end{proposition}

Let us for the moment assume
that we have proved Proposition \ref{final}.
Then it is not
difficult to derive Proposition \ref{pro-Galois-general}:

For any $t \in \mathbb{N}$, the cover $\mathcal{D}_t^0
\longrightarrow \mathcal{C}$ is a composite of the cover
$\mathcal{D}_t^0 \longrightarrow \mathcal{C}_{S_t^0}$ and the
cover $\mathcal{C}_{S_t^0} \longrightarrow \mathcal{C}$. By the
last item
of Proposition \ref{final} and the assumption that
$\pi_1(S)/\pi_1(S)^2$ is finite, there exists some $t_0 \in
\mathbb{N}$ such that for $t
> t_0$, the canonical morphism $\mathcal{D}_t^0 \longrightarrow
\mathcal{D}_{t_0}^0$ is an isomorphism.

Let $\mathcal{D}^0 := \mathcal{D}_{t_0}^0 =
\varprojlim
\mathcal{D}^0_{t}$, let $T:=S_{t_0}^0$. Then for all
$t \geq t_0$, $\mathcal{D}_t$ is a curve over $T$, and so is
$\mathcal{C}_T$. Let $\mathcal{D} :=
\varprojlim
\mathcal{D}_t$. Then $\mathcal{D}$ is a pro-Galois
curve cover of $\mathcal{C}_T$ of infinite degree. Moreover, its
Galois group is an extension of a group of the form
$(\mathbb{Z}/2\mathbb{Z})^r$ with $r \leq 4$ by
$\prod_{i=1}^{\infty}\Gal(\kappa(\mathcal{C}_i)|\kappa(\mathcal{C}_i^0))
\approx \prod_{i=1}^\infty A_{n_i}$. This implies Proposition
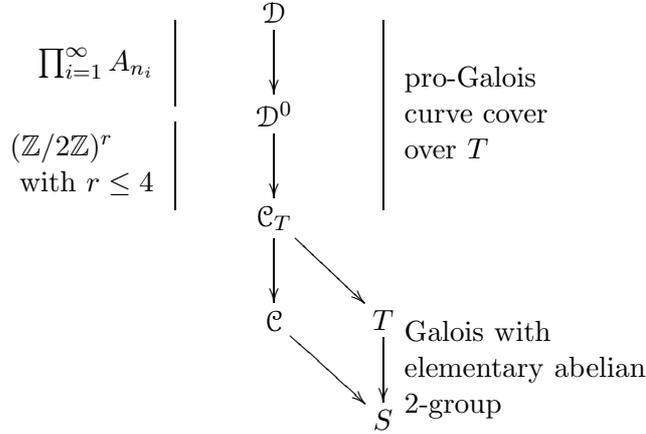
\ref{pro-Galois-general}.

\[ \xymatrix{
\ar@{-}[d]_{\begin{array}{l}{\prod_{i=1}^{\infty} A_{n_i}}\end{array}} & {\mathcal{D}} \ar[d] & \ar@{-}[dd]^{\begin{array}{l} \text{pro-Galois} \\ \text{curve cover} \\ \text{over } T \end{array}}  \\
\ar@{-}[d]_{\begin{array}{l}(\mathbb{Z}/2\mathbb{Z})^r \\ \text{ with } r \leq 4 \end{array}} &  {\mathcal{D}^0} \ar[d] \\
 & {\mathcal{C}_{T}} \ar[d] \ar[dr] &  \\
 & {\mathcal{C}} \ar[dr] & T \ar[d]\ar@{-}[d]^{\begin{array}{l} \text{Galois with} \\ \text{elementary abelian} \\ 2\text{-group} \end{array}}  \\
 &    & S
} \]

\smallskip

Now we give the proof of Proposition \ref{we-have-to-show}. It is
divided into several lemmata.

\subsubsection*{Proof of Assertion a)}

Let us recall Abhyankar's Lemma.

\begin{lemma}[Abhyankar's Lemma]
\label{usual-abhyankar}
Let $K$ be a field, $L|K,M|K$ finite separable extensions of $K$, $N=LM$ a compositum of $M$ and $L$
over $K$. Let $v$ be a discrete valuation of $N$, $v_M,v_L,v_K$ the restrictions of $v$ to $M,L,K$ respectively.
Assume that the extensions $v_M|v_K$ and $v_L|v_K$ are tame and that $e(v_M|v_K)|e(v_L|v_K)$.
Then $v|v_L$ is unramified.
\end{lemma}
For a \emph{proof} see \cite[Proposition III.8.9]{St} (the assumptions in \cite[Proposition III.8.9]{St} that $M|K$ be an extension of function fields in one variable and $v$ a valuation of function fields is not necessary).

\begin{lemma}
\label{finite,ramification...}
The morphisms $\widetilde{\mathcal{C}}_i \longrightarrow \mathcal{E}$ and
$\mathcal{D}_t \longrightarrow \mathcal{E}$ are finite, the branched points in codimension 1 are the generic
points of $\Delta$, and the corresponding ramification indices are 2.
\end{lemma}
\emph{Proof.}
The first statement follows from the following general fact: If one normalizes an integral noetherian scheme in a finite separable extension of its function field,
the canonical morphism is finite; cf.\ \cite[Chapter 4, Proposition 1.25]{Liu}. By Abhyankar's Lemma, the ramified points in codimension 1 are the generic points of $\Delta$ and the ramification indices divide 2. As the corresponding extensions of function fields are Galois, the ramification indices for points above the branched points in codimension 1 are 2.
\qed

\begin{lemma}
\label{etale-covers}
The morphisms $\mathcal{D}_t \longrightarrow \widetilde{\mathcal{C}}_0 = \mathcal{C}$ are \'etale covers.
\end{lemma}
\emph{Proof.} Let $t \in \mathbb{N}$ be fixed.

By Lemma \ref{finite,ramification...}, the branch loci of $\mathcal{D}_t \longrightarrow \mathcal{E}$
and $\mathcal{C} \longrightarrow \mathcal{E}$ are equal,
and the ramification indices of the branched points in codimension 1 are dividing 2.
As $L_0|\kappa(\mathcal{E})$ has degree 2,
this implies that for all points $x$ in codimension 1 of $\mathcal{E}$, all ramification indices of
$\mathcal{D}_t \longrightarrow \mathcal{E}$ at $x$ divide all ramification indices of
$\widetilde{\mathcal{C}}_0 \longrightarrow \mathcal{E}$ at $x$.

By Abhyankar's Lemma, $\mathcal{D}_t \longrightarrow \mathcal{C}$ is unramified at all points of codimension 1.
As $S$ is regular and $\mathcal{C}$ is smooth over $S$, $\mathcal{C}$ is regular (in particular normal).

By ``purity of the branch locus'' (\cite[Expos\'e X, 3.1]{SGA}), $\mathcal{D}_t \longrightarrow \mathcal{C}$
is unramified everywhere, and because $\mathcal{C}$ is normal, it is \'etale;
see \cite[Expos\'e I, Corollaire 9.11]{SGA}.
\qed

\begin{lemma}
\label{Galois-covers}
The morphisms $\mathcal{D}_t \longrightarrow \mathcal{C}$, $\mathcal{D}_t \longrightarrow \mathcal{D}_t^0$
and $\mathcal{D}_t^0 \longrightarrow \mathcal{C}$ are Galois covers.
\end{lemma}
\emph{Proof.}
By the last lemma, we already know that the morphisms $\mathcal{D}_t \longrightarrow \mathcal{C}$ are \'etale covers.
Moreover, the corresponding extensions of functions fields are Galois. As $\mathcal{C}$ (respectively $\mathcal{D}_t^0$)
is normal, this implies that the covers $\mathcal{D}_t \longrightarrow \mathcal{C}$ and
$\mathcal{D}_t \longrightarrow \mathcal{D}_t^0$ are Galois. This implies with \cite[Corollaire 3.4]{SGA}
that $\mathcal{D}_t^0 \longrightarrow \mathcal{C}$ is \'etale. It is Galois because the corresponding extension of
function fields is Galois (and $\mathcal{C}$ is normal).
\qed

\subsubsection*{Proof of Assertions b) and c)}

We need a group theoretical lemma.
\begin{lemma}
\label{groups-simple-surjective}
 Let $G_1, \dots, G_t$ be groups such that for any $i, j$ with $i \neq j$, there does not exist any non-trivial simple group which is a quotient group of both $G_i$ and $G_j$. Let $G$ be a group with surjective homomorphisms $p_i : G \longrightarrow G_i$. Then the induced homomorphism $p : G \longrightarrow \prod_{i=1}^t G_i$ is surjective. 
\end{lemma}
\emph{Proof.} The proof can be done by induction on $t$. For $t=2$, the proof is as follows.

Let $N:= \langle \ker(p_1) \cup \ker(p_2) \rangle$ -- this is a normal subgroup of $G$.
We have canonical surjective homomorphisms $G_i \simeq G/\ker(p_i) \longrightarrow G/N$.
By assumption, $G/N$ is trivial, that is, $G=N$.

This implies that $p_1|_{\ker(p_2)} : \ker(p_2) \longrightarrow G_1$ and $p_2|_{\ker(p_1)} :
\ker(p_1) \longrightarrow G_2$ are surjective. This in turn implies that the image of $p$
contains $G_1 \times \{ 1 \} $ and $\{ 1 \} \times G_2$, thus $p$ is surjective.
\qed

\medskip

Now we can proceed with the proof of Proposition \ref{final}.

By construction, $\kappa(\mathcal{D}^0_t)|\kappa(\mathcal{E})$ is
Galois with Galois group an elementary abelian 2-group. It follows
that the Galois group of $\mathcal{D}^0_t \longrightarrow
\mathcal{C}$ is an elementary abelian 2-group. This proves Assertion b).

The covers $\mathcal{D}_t \longrightarrow \mathcal{D}^0_t$ are more interesting:

\begin{lemma}
\label{A_n}
The restrictions $\Gal(\kappa(\mathcal{D}_t)|\kappa(\mathcal{D}^0_t)))
 \longrightarrow \Gal(\kappa(\widetilde{\mathcal{C}}_i)|\kappa(\widetilde{\mathcal{C}}_i^0))$ are surjective and induce an isomorphism
\[\Gal(\mathcal{D}_t \rightarrow \mathcal{D}^0_t) \simeq \Gal(\kappa(\mathcal{D}_t)|\kappa(\mathcal{D}^0_t)) \tilde{\longrightarrow} \prod_{i=1}^{t} \Gal(\kappa(\widetilde{\mathcal{C}}_i)|\kappa(\widetilde{\mathcal{C}}_i^0)) \approx \prod_{i=1}^t A_{n_i}.\]
\end{lemma}
\emph{Proof.} The group $\Gal(\kappa(\widetilde{\mathcal{C}}_0)|\kappa(\widetilde{\mathcal{C}}_0^0))$ is trivial, so let $i \geq 1$. The group
$\Gal(\kappa(\widetilde{\mathcal{C}}_i)|\kappa(\widetilde{\mathcal{C}}_i^0))$
is isomorphic to $A_{n_i}$, and the group
$\Gal(\kappa(\mathcal{D}_i^0)|\kappa(\widetilde{\mathcal{C}}_i^0))$
is an elementary abelian 2-group. As $n_i \geq 5$ (because $i \geq 1$), the
group $A_{n_i}$ is simple (see \cite[Theorem 3.15]{Ro}) and in
particular has no non-trivial elementary abelian 2-group as
quotient. By Lemma \ref{groups-simple-surjective}, the extensions
$\kappa(\mathcal{D}^0_t)$ and $\kappa(\widetilde{\mathcal{C}}_i)$
are linearly disjoint over $\kappa(\widetilde{\mathcal{C}}^0_i)$
(inside $\kappa(\mathcal{D}_t)$) (this will also be used in the proof of Lemma \ref{regular}). In particular, the restriction map
$\Gal(\kappa(\mathcal{D}^0_t) \,
\kappa(\widetilde{\mathcal{C}}_i)|\kappa(\mathcal{D}^0_t))
\longrightarrow
\Gal(\kappa(\widetilde{\mathcal{C}}_i)|\kappa(\widetilde{\mathcal{C}}_i^0))$
is an isomorphism.

Again by Lemma \ref{groups-simple-surjective}, the induced homomorphism
$\Gal(\kappa(\mathcal{D}_t)|\kappa(\mathcal{D}^0_t))$ \linebreak
$\longrightarrow \prod_{i=1}^t \Gal(\kappa(\widetilde{\mathcal{C}}_i)|\kappa(\widetilde{\mathcal{C}}_i^0))
\approx \prod_{i=1}^t A_{n_i}$ is surjective. It is obvious that it is injective.
\qed

\subsubsection*{Proof of Assertion d)}

Again we first need a group theoretical lemma.

\begin{lemma}
\label{groups-simple-injective}
Let $G_1, \ldots, G_t$ be finite groups such that any $i, j$ with $i \neq j$, there is no simple group which occurs as a composition factor of both $G_i$ and $G_j$. Let $G$ be a group with a homomorphism $\varphi : \prod_{i=1}^t G_i \longrightarrow G$ such that for $i = 1, \ldots, t$, the restriction of $\varphi$ to $G_i$ (regarded as a subgroup of $\prod_{i=1}^t G_i$) has kernel $N_i$. Then $\varphi$ has kernel $\prod_{i = 1}^t N_i$ (regarded as a subgroup of $\prod_{i=1}^t G_i$).
\end{lemma}
\emph{Proof.}
Recall that if $G$ is a finite group and $N \vartriangleleft G$ is a normal subgroup, then composition series of both $N$ and $G/N$ in a canonical way give rise to a composition series of $G$. In particular, the set of composition factors of $G$ is the union of the sets of composition factors of $N$ and $G/N$.

The assumption and this remark imply that for any $j =2, \ldots, t$, there is no simple group which occurs in the composition series of both $\prod_{i=i}^{j-1} G_i$ and $G_j$. Because of this, the general case follows by induction from the case for $t=2$.

In this case, the proof is as follows.

Obviously, $N_1 \times N_2$ in contained in the kernel of $\varphi$. The group $\varphi(G_1 \times \{1 \} ) \cap \varphi(\{1 \} \times G_2)$ is normal in both $\varphi(G_1 \times \{ 1 \})$ and $\varphi(\{ 1 \} \times G_2)$. 

However by assumption and the remark at the beginning of the proof, the groups $\varphi(G_1 \times \{ 1 \})$ and $\varphi(\{1 \} \times G_2)$ cannot contain a non-trivial common normal subgroup. This implies that the group $\varphi(G_1 \times \{1 \} ) \cap \varphi(\{1 \} \times G_2)$ is trivial. 

Now if $\varphi(g_1,g_2) = 1$, then $\varphi(g_1,1) = \varphi(1,g_2^{-1}) \in \varphi(G_1 \times \{1 \} ) \cap \varphi(\{1 \} \times G_2)$. Thus this has to be 1. By the definition of $N_1$ and $N_2$, $g_1 \in N_1$, $g_2 \in N_2$, that is, $(g_1,g_2) \in N_1 \times N_2$.
\qed

\begin{lemma}
\label{regular} Let $F_t$ be the algebraic closure of
$\kappa(S)$ in $\kappa(\mathcal{D}^0_t)$.
Then $\kappa(\mathcal{D}_t)|F_t$ is regular.
\end{lemma}
\emph{Proof.} 
As we have seen in the proof of Lemma \ref{A_n}, the fields
$\kappa(\mathcal{D}^0_t)$ and $\kappa(\widetilde{\mathcal{C}}_i)$
are linearly disjoint extensions of $\kappa(\widetilde{\mathcal{C}}^0_i)$. Let $F_{t,i}$ be the algebraic closure of $F_t$ in $\kappa(\mathcal{D}_t^0) \kappa(\widetilde{\mathcal{C}}_i)$. If we apply Lemma \ref{groups-simple-injective} to the restriction map $\text{Gal}(\kappa(\mathcal{D}_t^0) \kappa(\widetilde{\mathcal{C}}_i)|\kappa(\widetilde{\mathcal{C}}^0_i)) \longrightarrow \text{Gal}(F_{t,i} \kappa(\tilde{\mathcal{C}}^0_i)|\kappa(\tilde{\mathcal{C}}^0_i))$, we obtain that $F_{t,i} \kappa(\tilde{\mathcal{C}}^0_i) = F_t \kappa(\tilde{\mathcal{C}}_i^0)$. As $\kappa(\tilde{\mathcal{C}}^0_i)|\kappa(S)$ is regular, we obtain $F_{t,i} = F_t$, that is, $\kappa(\mathcal{D}_t^0) \kappa(\widetilde{\mathcal{C}}_i)|F_t$ is regular.

By the structure of
$\text{Gal}(\kappa(\mathcal{D}_t)|\kappa(\mathcal{D}_t^0))$ and
again Lemma \ref{groups-simple-injective}, this implies that
$\kappa(\mathcal{D}_t)|F_t$ is regular. \qed

\begin{lemma}
Let $q_t : \mathcal{D}^0_t \longrightarrow S$ be the structure morphism, and let $S_t^0$ be the normalization of $S$ in $F_t$. Then $S_t^0$ is equal
to $\mathbf{Spec}({q_t}_*(\mathcal{O}_{\mathcal{D}^0_t}))$, and the canonical
morphism $S_t^0 \longrightarrow S$ is an \'etale cover.
\end{lemma}
\emph{Proof.} As $\mathcal{D}^0_t$ is smooth and proper over $S$, for all $s
\in S$, the ring of global sections of the fiber $(\mathcal{D}_t^0)_s$
is a finite separable algebra over $\kappa(s)$. This implies with
\cite[III (7.8.7)]{EGA} that
$\textbf{Spec}({q_t}_*(\mathcal{O}_{\mathcal{D}^0_t}))
\longrightarrow S$ is an \'etale cover. Again by \cite[III
(7.8.7)]{EGA} applied to the generic point of $S$, one sees that
the total ring of fractions of
$\textbf{Spec}({q_t}_*(\mathcal{O}_{\mathcal{D}^0_t}))$ is $F_t$.
This implies that
$\textbf{Spec}({q_t}_*(\mathcal{O}_{\mathcal{D}^0_t}))$ is the
normalization of $S$ in $F_t$, that is, it is $S_t^0$; see \cite[Corollaire 10.2]{SGA}.
\qed

\medskip

We can consider $\mathcal{D}_t^0$ as an $S^0_t$-scheme (Stein
factorization); let $r_t : \mathcal{D}^0_t \longrightarrow S^0_t$
be the structure morphism. Then $\mathcal{O}_{S^0_t} =
{r_t}_*(\mathcal{O}_{\mathcal{D}^0_t})$, and $\mathcal{D}^0_t$ has
connected and non-empty geometric fibers over $S^0_t$; see
\cite[III (4.3.1), (4.3.4)]{EGA}.

By the universal property of the product, $\mathcal{D}^0_t
\longrightarrow \mathcal{C}$ factors through $\mathcal{C}_{S_t^0}
\longrightarrow \mathcal{C}$.

\begin{lemma}
$\mathcal{D}^0_t$ is a curve over $S_t^0$ and $\mathcal{D}^0_t \longrightarrow \mathcal{C}_{S_t^0}$ is an \'etale curve cover.
\end{lemma}
\emph{Proof.} The morphism $\mathcal{D}_t^0 \longrightarrow \mathcal{C}_{S_t^0}$ is an \'etale cover,
because $\mathcal{D}_t^0 \longrightarrow \mathcal{C}$ and $\mathcal{C}_{S_t^0}
\longrightarrow \mathcal{C}$ are \'etale covers. This implies that $\mathcal{D}^0_t$ is smooth and proper over $S_0^t$. It is a curve because it has connected and 1-dimensional geometric fibers. \qed

\medskip

As $S_t^0 \longrightarrow S$ is an \'etale cover and $F_t|\kappa(S)$ is Galois with
Galois group a quotient of that of $\kappa(\mathcal{C}_t^0)|\kappa(\mathcal{C})$
(as one easily sees) and $S$ is normal, $S_t^0 \longrightarrow S$ is a Galois cover
with Galois group an elementary abelian 2-group. The Galois group is canonically isomorphic
to that of $\mathcal{C}_{S^0_t} \longrightarrow \mathcal{C}$ because base-change does not
change the Galois group.

By Galois theory the cover $\mathcal{D}_t^0 \longrightarrow \mathcal{C}_{S_t^0}$
is Galois, and we have the exact sequence
\[ 1 \longrightarrow \Gal(\mathcal{D}_t^0 \rightarrow \mathcal{C}_{S^0_t})
 \longrightarrow \Gal(\mathcal{D}_t^0 \rightarrow \mathcal{C}) \longrightarrow
 \Gal(\mathcal{C}_{S_t^0} \rightarrow \mathcal{C}) \longrightarrow 1 \]
with $\Gal(S_t^0 \rightarrow S) \simeq \Gal({\mathcal{C}_{S_t^0}
\rightarrow \mathcal{C}})$, where all groups are elementary
abelian 2-groups. By pull-back to a geometric fiber, we see that
$\Gal(\mathcal{D}_t^0 \rightarrow \mathcal{C}_{S^0_t})$ is a group
of the form $(\mathbb{Z}/2\mathbb{Z})^r$ with $r \leq 2 \cdot g_{\mathcal{C}} = 4$.

Now let $S_t$ be analogously defined to $S_t^0$ for
$\mathcal{D}_t$. Then $\mathcal{D}_t$ has connected geometric
fibers over $S_t$ and $S_t$ is \'etale over $S_t^0$. By Lemma
\ref{regular} and \cite[III (7.8.7)]{EGA}, the generic fiber of
$\mathcal{D}_t$ over $S_t^0$ is geometrically connected. On the
other hand, again by \cite[III (7.8.7)]{EGA} the number of
connected components in each geometric fiber of $\mathcal{D}_t$
over $S^0_t$ is equal to the degree of $S_t$ over $S_t^0$. This
implies that $S_t \tilde{\longrightarrow} S_t^0$. The fact that
$\mathcal{D}_t \longrightarrow \mathcal{D}^0_t$ is an \'etale cover
(Lemma \ref{Galois-covers}) implies that $\mathcal{D}_t$ is a curve over $S_t^0$. This completes the proof of Proposition~\ref{final}.

\section{Complementary Results}

\subsection{Abelian Surfaces Isogenous and Isomorphic to a Product of Elliptic Curves}
\label{isogenous-isomorphic}

Recall that if $\mathcal{C} \longrightarrow \mathcal{E}$ is a minimal cover of degree $n$ over some
scheme $S$, the Jacobian of $\mathcal{C}$ is as principally polarized
abelian variety isomorphic to $((\mathcal{E} \times \mathcal{E}')/\text{Graph}(-\psi),
\lambda)$ where $\mathcal{E}'$ is an elliptic curve over $S$, $\psi : \mathcal{E}[n]
\longrightarrow \mathcal{E}'[n]$ is an isomorphism which is anti-isometric with respect to
the Weil pairing and $\lambda$ is the (principal) polarization whose pull-back to $\mathcal{E}
 \times \mathcal{E}'$ is $n$ times the product polarization; see \cite{KaHur}.

An important special case is that $\psi$ is induced by an isogeny
$\tau : \mathcal{E} \longrightarrow \mathcal{E}'$ (necessarily of degree coprime to $n$).
 In this case $(\mathcal{E} \times \mathcal{E}')/\text{Graph}(-\psi)$
 is in fact always isomorphic to $\mathcal{E} \times \mathcal{E}'$ as can
 be seen by the following easy generalization of the proof of Proposition \ref{isomorphic}.
\smallskip

Let $n$ is some natural number, let $\mathcal{E}$, $\mathcal{E}'$ two elliptic curves over a scheme $S$,
let $\tau : \mathcal{E} \longrightarrow \mathcal{E}'$ be any isogeny of degree $N$ coprime to $n$. (In Proposition \ref{isomorphic} we treated the case that $n$ is 2 and $N$ is odd.) As in Section \ref{over Y_0(N)}, let $\psi := \tau|_{\mathcal{E}[n]}$.

Then the isogeny $\Phi : \mathcal{E} \times \mathcal{E}' \longrightarrow \mathcal{E} \times \mathcal{E}'$
given by the matrix $\left( \begin{array}{cc} n & 0 \\ \tau & 1 \end{array} \right)$ has kernel $\Graph(-\psi)$,
thus it induces in isomorphism
\begin{equation}
\label{isog-isom-general}
(\mathcal{E} \times \mathcal{E}')/\Graph(-\psi) \tilde{\longrightarrow} \mathcal{E} \times \mathcal{E}'.
\end{equation}
Assume now that $\psi$ is anti-isometric with respect to the Weil pairing, i.e.\ that $N$ is congruent to
$-1$ modulo $n$. Then $(\mathcal{E} \times \mathcal{E}')/\text{Graph}(-\psi)$ has a principal polarization
$\lambda$ whose pull-back to $\mathcal{E} \times \mathcal{E}'$ is $n \, \id_{\mathcal{E} \times \mathcal{E}'}$;
see \cite[Proposition 5.7]{KaHur}. (As in Section \ref{over Y_0(N)}, we identify elliptic curves (and products
of elliptic curves) with
their duals, so that the canonical (product) polarizations become
the identity.)

Under isomorphism (\ref{isog-isom-general}), the polarization $\lambda$ corresponds to the principal
polarization $\widetilde{\lambda}$ on $\mathcal{E} \times \mathcal{E}'$ whose pull-back with $\Phi$ is
also $n \, \id_{\mathcal{E} \times \mathcal{E}'}$, i.e.\ we have
$\widehat{\Phi} \, \widetilde{\lambda} \, \Phi = n \, \id_{\mathcal{E} \times \mathcal{E}'}$. It follows that $\widetilde{\lambda} = n \, \widehat{\Phi}^{-1} \, \Phi^{-1}$,
and consequently $\widetilde{\lambda}$ is given by the matrix
\begin{equation}
\label{induced-polarization}
\frac{1}{n} \left( \begin{array}{cc} 1 & - \widehat{\tau} \\ 0 & n \end{array} \right) \,
\left( \begin{array}{cc} 1 & 0 \\ - \tau & n  \end{array} \right) = \left( \begin{array}{cc}
\frac{1 + N}{n} & - \widehat{\tau} \\ - \tau & n \end{array} \right).
\end{equation}

In particular, any principally polarized abelian surface considered in \cite[Proposition 3.1]{FrHong} is
isomorphic to a product of elliptic curves with a polarization given by (\ref{induced-polarization}).

\subsection{The Number Field Case}
\label{number-case}
It is interesting to note that the assumptions of Proposition \ref{pro-Galois-general} cannot be satisfied for
$S$ equal to on open part of $\Spec(\mathcal{O}_K[1/2])$, where $\mathcal{O}_K$ is the principal order
in a number field $K$. Indeed, this is contradicted by the following proposition.
\begin{proposition}
Let $K$ be a number field, let $S$ an open part of $\Spec(\mathcal{O}_K[1/2])$, and let
$\mathcal{E}$ be an elliptic curve over $S$. Then there does not
exist a sequence $(\pi_i)_{i \in \mathbb{N}_0}$ of minimal covers
$\pi_i : \mathcal{C}_i \longrightarrow \mathcal{E}$ (where the
$\mathcal{C}_i$ are curves of genus $2$) with pairwise distinct
degrees and the same branch locus which
is \'etale of degree $2$ over $S$.
\end{proposition}
\emph{Proof.} By Faltings' proof of the Shafarevich Conjecture (\cite{Fal}),
there exist only finitely many isomorphism classes of curves of
genus $2$ over $S$. It thus suffices to prove the proposition
under the assumption that the $\mathcal{C}_i$ are equal to each
other. Now the conclusion is implied by the following proposition.
\qed

\begin{proposition}
\label{not-equal-number-field}
Let $K$ be a number field, $E$ an elliptic curve and
$C$ a curve of genus 2 over $K$. Then there does not exist a sequence $(\pi_i)_{i \in \mathbb{N}_0}$ of
minimal covers $\pi_i : C \longrightarrow E$ with pairwise distinct degrees and the same
branch locus which has degree 2 over $K$.
\end{proposition}
\emph{Proof.}
Assume there exists a sequence $(\pi_i)_{i \in \mathbb{N}_0}$ of minimal covers with pairwise distinct degrees.
Let $\Delta$ be the common branch locus, let $V_i$ the ramification locus of $\pi_i$. By the assumption and
the Hurwitz genus formula both $\Delta$ and $V_i$ have degree 2, thus the canonical maps
$V_i \longrightarrow \Delta$ are isomorphisms. As by Proposition \ref{ramification-alternative}
the $V_i$ are pairwise disjoint, $C$ has infinitely many $\Delta$-valued points. This contradicts that by Faltings' proof of the Mordell Conjecture (\cite{Fal}), any curve of genus $\geq 2$ over a number field has only finitely many rational points.\qed

\medskip

One might ask whether there exists an elliptic curve $E$ over a number field $K$ and a sequence
 $(\pi_i)_{i \in \mathbb{N}}$ of minimal covers $\pi_i : C_i \longrightarrow E$ (where the $C_i$ are
 curves of genus $2$ over $K$)
 with pairwise distinct degrees and equal branch loci. With similar (but easier) arguments as in
  Section \ref{infinite-degree}, the existence of such a sequence would lead to a curve of genus 2 over the maximal
  elementary abelian 2-extension of $K$ which has a pro-Galois curve cover of infinite degree.

There is the following conjecture closely related to the height conjecture for
elliptic curves and so to the ABC conjecture.

\paragraph{Conjecture}
Let $E$ be an elliptic curve over a number field $K$. Then there is a number
$n_0(K,E)$ such that for all elliptic curves $E'$
over $K$ with $E[n]$ isomorphic to $E'[n]$ for some $n>n_0(K,E)$ it
follows that $E$ and $E'$ are isogenous (over $K$).

\medskip

This conjecture is equivalent to Conjecture 5 in \cite{FrBos} in the special case that the base field is a number field. (We use that the Faltings height over a number field is effective and that given two elliptic curves $E$ and $E'$ over a number field such that for infinitely many natural numbers $n$, $E[n]$ is isomorphic to $E'[n]$, $E$ and $E'$ are isogenous.)

\begin{proposition}
Let $K$ be a number field, $E$ an elliptic curve over $K$. Then \emph{under the assumption of the above conjecture},
 there are only finitely many isomorphism classes of curves of genus $2$ occurring as minimal covers of $E$ with covering degree $> n_0(K,E)$.
\end{proposition}
\emph{Proof.} 
Let $C \longrightarrow E$ be a genus 2 cover with covering degree $n > n_0(K,E)$, $J_C \simeq (E \times E')/\Graph(-\psi)$, say. Then we have the isomorphism $\psi: E[n] \simeq E'[n]$, thus by the conjecture $E$ is isogenous to $E'$. This means that $J_C$ is isogenous to $E \times E$. By Faltings' results, there are only finitely many abelian surfaces isogenous to $E \times E$, and any of these finitely many abelian surfaces has up to isomorphism only finitely many principal polarizations. The result now follows by Torelli's Theorem.
\qed

\medskip

Together with Proposition \ref{not-equal-number-field}, this proposition implies:

\begin{proposition}
Let $K$ be a number field, $E$ an elliptic curve over $K$. Then \emph{under the assumption of the above conjecture},
 there does not exist a sequence $(\pi_i)_{i \in \mathbb{N}}$ of minimal covers
  $\pi_i : C_i \longrightarrow E$ (where the $C_i$ are curves of genus $2$) with pairwise distinct degrees and the
  same branch locus which has degree 2 over $K$.
\end{proposition}

\bibliographystyle{plain}
\bibliography{genus-2}

\begin{thebibliography}{10}

\bibitem{Di}
C.~Diem.
\newblock {Families of elliptic curves with genus 2 covers of degree 2}.
\newblock {\em Collect.\ Math.}, 57:1--25, 2006.

\bibitem{Fal}
G.~Faltings.
\newblock {Endlichkeitss\"atze f\"ur abelsche Varietäten \"uber Zahlkörpern}.
\newblock {\em Invent.\ Math.}, pages 349--366, 1983.

\bibitem{FC}
G.~Faltings and C.-L. Chai.
\newblock {\em {Degeneration of Abelian Varieties}}.
\newblock Springer-Verlag, Berlin, 1980.

\bibitem{FrHong}
G.~Frey.
\newblock On elliptic curves with isomorphic torsion structures and
  corresponding curves of genus {$2$}.
\newblock In {\em Elliptic curves, modular forms, \& Fermat's last theorem
  (Hong Kong, 1993)}, pages 79--98. Internat. Press, Cambridge, MA, 1995.

\bibitem{FrBos}
G.~Frey.
\newblock {On Ternary Equations of Fermat Type and Relations with Elliptic
  Curves}.
\newblock In G.~Cornell, J~Silverman, and G.~Stevens, editors, {\em {Modular
  Forms and Fermat's Last Theorem}}, pages 527--548. Springer-Verlag, 1997.

\bibitem{FK-proj}
G.~Frey and E.~Kani.
\newblock Projective $p$-adic representations of the $k$-rational geometric
  fundamental group.
\newblock {\em Arch.\ Math.}, 77(32-46), 2001.

\bibitem{FKV}
G.~Frey, E.~Kani, and H.~V{\"o}lklein.
\newblock Curves with infinite {$K$}-rational geometric fundamental group.
\newblock In {\em Aspects of Galois theory (Gainesville, FL, 1996)}, volume 256
  of {\em London Math. Soc. Lecture Note Ser.}, pages 85--118. Cambridge Univ.
  Press, Cambridge, 1999.

\bibitem{SGA}
A.~Grothendieck.
\newblock {\em Rev\^etements \'etales et groupe fondamental (SGA I)}, volume
  1960/61 of {\em S\'eminaire de G\'eom\'etrie Alg\'ebrique}.
\newblock Institut des Hautes \'Etudes Scientifiques, Paris.

\bibitem{EGA}
A.~{Grothendieck with J. Dieudonn\'e}.
\newblock {El\'ements de G\'eom\'etrie Al\-g\'e\-brique (I-IV)}.
\newblock {ch.\ I: Springer-Verlag, Berlin, 1971; ch. II-IV: Publ.\ Math.\
  Inst.\ Hautes Etud.\ Sci. 8,11, 17, 20,24, 28, 32, 1961-68}.

\bibitem{Ha}
R.~Hartshorne.
\newblock {\em Algebraic Geometry}.
\newblock Springer, New York, 1977.

\bibitem{Ih}
Y.~Ihara.
\newblock On unramified extensions of function fields over finite fields.
\newblock In {\em {Galois Groups and their representations}}, volume~2 of {\em
  Adv.\ Study Pure Math.}, pages 89--97. 1983.

\bibitem{KaHur}
E.~Kani.
\newblock Hurwitz spaces of genus 2 covers of an elliptic curve.
\newblock {\em Collect.\ Math.}, 54(1):1--51, 2003.

\bibitem{Lang-fundamentals}
S.~Lang.
\newblock {\em {Fundamentals of Diophantine Geometry}}.
\newblock Springer-Verlag, 1983.

\bibitem{La}
S.\ Lang.
\newblock {\em {Algebra (Third Edition)}}.
\newblock Addison-Wesley, 1993.

\bibitem{Liu}
Q.~Liu.
\newblock {\em {Algebraic Geometry and Arithmetic Curves}}.
\newblock Oxford University Press, Oxford, 2001.

\bibitem{Me}
L.~Merel.
\newblock Bornes pour la torsion des courbes elliptiques sur les corps de
  nombres.
\newblock {\em Invent.\ Math.}, 124:437--449, 1996.

\bibitem{Mi-EC}
J.~Milne.
\newblock {\em Etale Cohomology}.
\newblock Princeton University Press, 1980.

\bibitem{Mi-JV}
J.~Milne.
\newblock Jacobian varieties.
\newblock In G.~Cornell and J.~Silverman, editors, {\em Arithmetic Geometry},
  pages 167--212. Springer-Verlag, 1986.

\bibitem{Mu-GIT}
D.~Mumford.
\newblock {\em Geometric Invariant Theory}.
\newblock Springer-Verlag, Berlin, 1965.

\bibitem{Mu-AV}
D.~Mumford.
\newblock {\em {Abelian Varieties}}.
\newblock Tata Institute for Fundamental Research, 1970.

\bibitem{Ro}
J.~Rotman.
\newblock {\em An introduction to the theory of groups}.
\newblock Wm.\ C.\ Brown, 1994.

\bibitem{St}
H.~Stichtenoth.
\newblock {\em Algebraic Function Fields and Codes}.
\newblock Springer-Verlag, Berlin, 1993.

\bibitem{Wi}
H.~Wielandt.
\newblock {\em Finite Permutation Groups}.
\newblock Academic Press, New York, 1964.

\end{thebibliography}

\vspace{4 ex}
\noindent
Claus Diem. Universit\"at Leipzig, Fakult\"at f\"ur Mathematik und Informatik, Johannisgasse 26, 04109 Leipzig, Germany. email: diem@math.uni-leipzig.de\\[1 ex]
Gerhard Frey. Universit\"at Duisburg-Essen, Institut f\"ur Experimentelle Mathematik, Ellernstr. 29, 45326 Essen, Germany. email: frey@iem.uni-due.de

\end{document}